\documentclass[opre,nonblindrev]{mystyle} % current default for manuscript submission
\usepackage{amsmath}
\usepackage[ampersand]{easylist}
\usepackage{mathtools}
\usepackage{accents}
\usepackage{tabu}
\usepackage{amssymb}
\usepackage{bbm}
\usepackage{dsfont}
\usepackage{bbold}
\usepackage[numbers]{natbib}
\usepackage{multicol}
\usepackage{enumerate}
\usepackage{hyperref}
\usepackage{pbox}
\usepackage[inline]{enumitem}
\usepackage{algorithm}
\usepackage{units} %for \nicefrac
\usepackage[normalem]{ulem}

\usepackage{natbib}
 \bibpunct[, ]{(}{)}{,}{a}{}{,}%
 \def\bibfont{\small}%
 \def\bibsep{\smallskipamount}%
 \def\bibhang{24pt}%
\TheoremsNumberedThrough     % Preferred (Theorem 1, Lemma 1, Theorem 2)
\ECRepeatTheorems

\EquationsNumberedThrough    % Default: (1), (2), ...
\usepackage{color}
\definecolor{turquoise}{rgb}{0.19, 0.84, 0.78}
\newcommand{\blue}{\textcolor{black}}  %% just replace blue with black !

\newcommand{\sbf}{\boldsymbol}
\newcommand{\bs}{\boldsymbol}
\newcommand{\mb}{\mathbf}
\newcommand{\dd}{decision-dependent }
\newcommand{\ra}{\rightarrow}
\newcommand{\mU}{\mathcal{U}}

\usepackage{booktabs}
\usepackage{tabulary}
\usepackage{float}
\newcolumntype{K}[1]{>{\centering\arraybackslash}p{#1}}
\makeatletter
\newcommand*\bigcdot{\mathpalette\bigcdot@{.5}}
\newcommand*\bigcdot@[2]{\mathbin{\vcenter{\hbox{\scalebox{#2}{$\m@th#1\bullet$}}}}}
\makeatother

\begin{document}
%%%%%%%%%%%%%%%%

%\setcounter{page}{1}
%\DOI{10.1287/xxxx.0000.0000}%

\RUNAUTHOR{Nohadani and Sharma}

\RUNTITLE{Optimization under Decision-Dependent Uncertainty}

\TITLE{Optimization under \\ Decision-Dependent Uncertainty}

\ARTICLEAUTHORS{%
\AUTHOR{Omid Nohadani, Kartikey Sharma}
\AFF{Department of Industrial Engineering and Management Sciences, Northwestern University, Evanston, IL 60208, \EMAIL{\href{nohadani@northwestern.edu}{nohadani@northwestern.edu}, \href{kartikeysharma2014@u.northwestern.edu}{kartikeysharma2014@u.northwestern.edu}}} 
%\EMAIL{nohadani@northwestern.edu, kartikeysharma2014@u.northwestern.edu}} 
} % end of the block

% REQUIRED
\ABSTRACT{%
The efficacy of robust optimization spans a variety of settings with uncertainties bounded in predetermined sets.
In many applications, uncertainties are affected by decisions and cannot be modeled with current frameworks.
This paper takes a step towards generalizing robust linear optimization to problems with \dd uncertainties.
In general settings, we show these problems to be NP-complete.
To alleviate the computational inefficiencies, we introduce a class of uncertainty sets whose size depends on binary decisions.
We propose reformulations that improve upon alternative standard linearization techniques.
To illustrate the advantages of this framework, a shortest path problem is discussed, where the uncertain arc lengths are affected by decisions. 
Beyond the modeling and performance advantages, the proposed notion of proactive uncertainty control also mitigates over conservatism of current robust optimization approaches.
}

\HISTORY{revised January 3, 2018} %This paper has not been submitted yet.}

\maketitle

\section{Introduction}
\label{sec:intro}

The two well-established approaches of optimization under uncertainty are stochastic and robust optimization. 
Stochastic optimization (SO) can be used when the distribution of the uncertainty is available~\citep{Shapiro}.
When uncertainties can be regarded as residing in sets, robust optimization
(RO) is a computationally attractive alternative~\citep{Ben-Tal2009, Bertsimas2011}.
The method of RO has been extended considerably and applied to problems ranging from portfolio management~\citep{ghaoui2003worst}, to healthcare~\citep{chu2005robust}, to electricity systems~\citep{lorca2014multistage}, and to engineering design~\citep{bertsimas2010robust}. 

RO employs uncertainty sets that are predetermined and, hence, \emph{exogenous}.
For instance, temporal changes to the uncertainty can be explicitly modeled via time-dependent uncertainty sets~\citep{nohadani2017robust}.
In many real-world problems, however, the uncertainty can be affected by decisions.
In such \dd cases, the uncertainty set is \emph{endogenous}.
Despite the wide prevalence of such uncertainties in real-world settings, these problems have not received much attention in the past, largely due to computational intractabilities.
In this paper, we take a first step towards robust linear optimization problems with endogenous uncertainties and provide a class of uncertainty sets, whose reformulations improve over standard techniques.
Specifically, we study a single-stage RO problem with \dd uncertainty sets
\begin{equation}
\tag{RO-DDU}
\begin{aligned} 
\label{prob:ro_ddu}
\min_{\mb{x},\mb{y}} \>&\> \mb{c}^\top \mathbf{x} + \mb{f}^\top \mathbf{y}\\
\text{s.t.} \>&\> \mb{a}_{i}^\top \mathbf{x} + \bs{\xi}_i^\top \mathbf{y} \leq b_{i} \quad \forall \bs{\xi}_{i} \in \mathcal{U}_{i}(\mathbf{x}) \subseteq \mathbb{R}^{n} \quad \forall i = 1,\dots, m,
\end{aligned}
\end{equation}
where \(\mathbf{x}\in\mathbb{R}^n\) and \(\mathbf{y}\in\mathbb{R}^n\) represent decision
variables, which need to satisfy each constraint \({i=1, \dots, m}\) for every
realization from the set \(\mathcal{U}_{i}(\mathbf{x})\) that bounds the uncertain
parameter \(\bs{\xi}_i\).
Further, the parameters defining \(\mathcal{U}_{i}(\mathbf{x})\) depend on decisions  \(\mathbf{x}\).
%% content
We first study the complexity of~\eqref{prob:ro_ddu} for polyhedral \(\mathcal{U}_{i}(\mathbf{x})\).
We then assume $\mb{x}$ is binary and provide reformulations for special classes of polyhedral and conic uncertainty sets and conclude with numerical experiments.

To show the range of applicability of this model, we illustrate two examples.

\noindent\emph{Example 1: Uncertainty Reduction. $\quad$}
\blue{
In facility location or inventory management problems with uncertain demand, the uncertainty can be reduced by spending resources to acquire information.
Similarly, in healthcare problems, additional medical tests can improve the diagnosis.
This type of uncertainty reduction is characteristic of many real-world problems.
In order to improve solutions, decisions on uncertainty reduction have to be included into the optimization problem, making the uncertainty a function of decisions on acquiring additional information. 
}

\begin{figure}[h!]
  \begin{minipage}{0.45\textwidth}
      \includegraphics[width=1\linewidth]{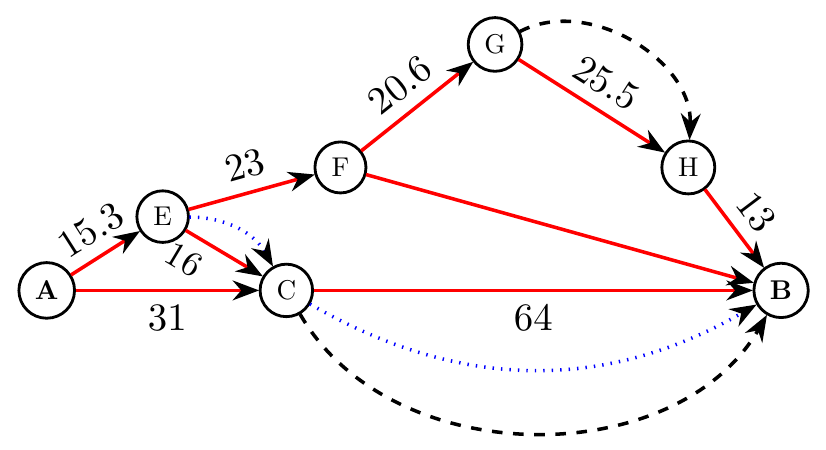}
  \end{minipage}
\hspace{2mm}
\begin{minipage}{0.52\textwidth}
\resizebox{\columnwidth}{!}{%
\begin{tabular}{|c|c|c|l|}
\hline
\textbf{Shortest Path} & \textbf{Path} & \textbf{Nominal}                    & \textbf{Worstcase}                        \\ \hline
Nominal                        & A$-$C$-$B         & $ 95$                          & $31 + 1.5\times64 =$\\&&& $\textbf{127}$                         \\ \hline
Robust                         & A$-$E$-$F$-$G$-$H$-$B   & \(97.4\) & $15.3 + 23 + 20.6 +$\\&&& $1.5\times25.5 + 13$\\&&& $= \textbf{110.15}$ \\ \hline
Endogenous          & A$-$E$-$C$-$B       & $95.3$             & $15.3 + 1.4\times16 +$\\  Robust &&& $1.1\times64 = \textbf{108.1}$           \\ \hline
\end{tabular}%
}

\end{minipage}
\caption{\label{fig:ex1}Shortest path on a network. Nominal lengths are labeled. Worst-case and reduced-case lengths are displayed with dashed and dotted lines. The table shows the lengths in different settings.}
\end{figure}

\noindent\emph{Example 2: Shortest Path on a Network. $\quad$}
Consider the graph in Figure~\ref{fig:ex1} with the arcset \(\mathcal{A}\) and let the uncertain length for any arc \(e\) be \(d_e = \bar{d}_e(1 + 0.5\xi_e)\), where $\bar{d}_e$ denotes the nominal value.
The uncertain vector \(\bs{\xi}\) lies in the uncertainty set \(\mU(\mb{x}) = \left\{\bs{\xi} \mid 0 \leq \xi_e \leq 1- 0.8x_e \;\;\forall e,\; \sum_{e \in \mathcal{A}}\xi_e \leq 1\right\}\).
The binary decision $x_e$ determines whether to reduce the maximum possible uncertainty \(\xi_e\) to 0.2 (\(x_e = 1\)) or leave it at 1 (\(x_e = 0\)).
For simplicity, we assume the reduction to be possible for at most one of the arcs.

Figure~\ref{fig:ex1}  displays a network with source node A and destination B. 
In the worst-case, the nominally shortest path lengthens to 127 units.
RO optimizes against this case, improving the worst-case length. 
If it is permitted to reduce the uncertainty of an arc, then A$-$E$-$C$-$B is selected with $x_{C-B}=1$ and the worst-case path becomes 108.5.
This example demonstrates that \dd sets can be leveraged to model decisions that mitigate the
worst-case scenario.

The contributions of this paper can be summarized as follows:
\begin{enumerate}
\item We study robust linear optimization problems with a polyhedral decision\-dependent set for the uncertain parameters.
We prove such problems to be NP-complete.
We also show that when decisions that influence the uncertainties are binary, the problem can be reformulated as a mixed integer optimization problem. 

\item For binary $\mb{x}$, we provide a class of uncertainty sets for which a more efficient reformulation of the \dd RO problem is possible.
The set structure and the nature of decision dependence are leveraged to provide reformulations with fewer constraints. 

\item We provide an improvement to \mbox{Big-M} linearization for bilinear terms which can reduce the number of constraints.

\end{enumerate}
This work also showcases the advantages that can be gained in both stochastic and robust optimization by proactively controlling uncertainties. 

We also emphasize what this paper fails to address.
Reformulations for continuous decisions influencing the uncertainty are not provided.
Furthermore, the primary problem in this paper is a static optimization problem,
i.e., the decisions do not adapt to uncertainty realizations. 
In fact, it is the uncertainty set and the corresponding worst-case realization
that are affected by decisions. 

Section~\ref{sec:gen_ddu} discusses the complexity of the \dd  robust linear optimization problem. 
Section~\ref{sec:red_ddu} introduces a class of uncertainty sets which allow improved reformulations. 
Section~\ref{sec:ext_gen} provides a comparison to the corresponding \mbox{Big-M} formulation.
It also provides methods to improve these standard techniques.
A numerical experiment is discussed in Section~\ref{sec:num_exp} to illustrate the advantages of the \dd setting and to computationally compare the three formulations.

\emph{Notation. $\quad$}
Throughout this paper, we use bold lower  and uppercase letters to denote vectors and matrices.
Scalars are marked in regular font.
All vectors are column vectors and the vector of ones is denoted by $\mb{e}$.
Furthermore, $diag({\bigcdot})$ denotes a diagonal matrix with $\bigcdot$ on the diagonal and zeros elsewhere.
For any given matrix \(\mb{A}\), the \(i^{th}\) row is denoted by \(\mb{A}_{i,\bigcdot}\)
and the \(j^{th}\)  column is denoted by \(\mb{A}_{\bigcdot,j}\).
The problems have \(m\) constraints indexed by \(i\). LHS denotes left-hand-side and RHS denotes right-hand-side.
We use the phrases ``decision-dependent'' and ``endogenous'' interchangeably.
Similarly, we refer to variables affecting an uncertainty set as influence variables.

%%%%%%%%%%%%%%%%
\section{Background}
\label{sec:background}

In the following, we first review endogenous settings in SO before discussing RO approaches.

The notion of endogenous uncertainty in SO generally corresponds to scenario trees, where decisions determine the probabilities.
For example, \citet{Jonsbraten1998} consider the cost of an item to remain uncertain until it is produced. 
The probability distribution depends upon which item is to be produced and when.
\citet{Goel2004} address the problem of offshore oil and gas planning, with the objective of maximizing revenues and investments over a period of time, when the recovery and size of oil fields are not known in advance. 
They provide a disjunctive formulation that is solved by a decomposition algorithm.
This approach is extended to a multistage SO problem for optimal production scheduling, that minimizes cost while satisfying the demand for different goods~\citep{Goel2006}.
For package sorting centers, \citet{novoa_sorting_centers} seek to balance the flow across working stations.
Capacities are modeled via budgeted uncertainties where the budget is a function of workstation allocation.
These and other approaches address endogenous uncertainties probabilistically.

In RO, the endogenous nature of uncertainty is imposed directly on the uncertainty set
itself.
For example, \citet{Spacey2012} address a software partitioning problem, where code segments are assigned to different computing nodes to reduce runtime with uncertain execution order and for unknown frequency of segment calls.
They employ tailored \dd uncertainty sets.
Such sets also occur as a result of reformulations.
For example, \citet{Hanasusanto2015} use a finite adaptability approximation to adjustable robust optimization (ARO),
as introduced by \citet{Bertsimas2010}, and consider optimization problems with binary recourse decisions.
For problems with uncertain objective and constraints, they provide a formulation with \dd uncertainty sets before finally reformulating it as a MILP.
\citet{poss2013robust,poss2014robust} considers combinatorial optimization problems with budgeted uncertainty sets.
This extends the work of~\citet{Bertsimas2004} to \dd budgets.
These works focus on budget uncertainty sets with limited discussion on general sets.
On the other hand, for a dynamic pricing problem with learning,~\citet{Bertsimas} consider 1 or \(\infty\)-norm uncertainty sets for price-dependent demand.
Specifically, the uncertain demand curve is driven by past realizations of price-demand pairs.
Since the price is a decision variable, this leads to \dd uncertainty sets.
In the context of robust scheduling problems, \citet{vujanic2016robust} consider a \dd uncertainty set which is a vector sum of a collection of sets.
The sets in the vector combination are selected by a decision which is a part of the original problem.
They probe the performance of an affine policy for the problem. 
More recently, \dd sets were studied in the context of control problems with \emph{primitive} uncertainty sets~\cite{ZhaEtal:2017:IFA_5415}.
Note that in all approaches to date, the decision dependence is modeled in a specific context, often driven by an application.

%% \section{conservatism}
The journey of RO has also included measures to reduce conservatism.
The original RO formulation by~\citet{Soyster1973} produced over conservative solutions for many applications due to the use of box uncertainties. 
Later, \citet{Ben-Tal1999} provided less conservative solutions by using general polyhedral and ellipsoid uncertainty sets.
ARO models~\cite{Ben-Tal2004} and decision rule approximations took another step in this direction by allowing decisions to depend on the realizations~\citep{Iancu2010,georghiou2015generalized}.
In this vein, \dd uncertainty sets offer a new avenue to reduce the level of conservatism. 
For example, \citet{poss2013robust} decreases it for cardinality constrained sets.
This work also motivates the notion of \emph{proactive uncertainty control} by using \dd sets to enable deliberate uncertainty reduction.

%%%%%%%%%%%%%%%%%

%%%%%%%%%%%%%%%%%%%%%%%%%%
\section{General Decision Dependence}
\label{sec:gen_ddu}
Robust linear optimization problems encompass a wide variety of applications,
in portfolio optimization, healthcare, inventory management, and routing, amongst others.
The tractability of robust linear programs provides a suitable starting point to analyze the
complexity of RO problems with \dd uncertainty.
Here, we investigate a robust linear optimization problem as in~\eqref{prob:ro_ddu}.
The underlying uncertainty set is endogenous and defined as follows.
\begin{definition}
  The set with constraint matrix \(\mb{D}\), constant vector \(\mb{d}\), and decision coefficient matrix \(\bs{\Delta}\) given by
  \begin{equation*}
  \label{U-P}
  \mathcal{U}^P(\mathbf{x}) = \{\bs{\xi} \>|\> \mb{D}\bs{\xi} \leq \mb{d} + \bs{\Delta}\mathbf{x}\}
  \end{equation*}
  is a \textbf{polyhedral} uncertainty set with affine decision dependence.
\end{definition}
Note that $\bs{\Delta}$ determines the influence of $\mb{x}$ on the set and can be estimated from the data or from the context of an application.
In Section~\ref{sec:num_exp}, we quantify it for a specific application.

The following theorem shows that RO problems with \dd sets cannot be reformulated in a tractable fashion, a departure from standard RO problems.
This occurs despite the fact that linear programs with polyhedral uncertainty sets have tractable robust counterparts.

\label{sec:NPC}
\begin{theorem}
\label{prop_npc}
The robust linear problem~\eqref{prob:ro_ddu} with uncertainty set $\mU^{P}$ is NP-complete.
\end{theorem}
\begin{proof}
The proof follows the following steps:
\begin{enumerate}
\item Consider an instance of the 3-Satisfiability problem (3-SAT) for a set of literals \(N = \{1,2, \dots, n\}\) and \(m\) clauses, which seeks to find a solution \(\mb{x} \in \{0,1\}^{n}\) that satisfies 
  \[x_{i_1} + x_{i_2} + (1-x_{i_3}) \geq 1 \;\; \text{for } m \text{ clauses and } i_1, i_2, i_3\in\{1,\dots,n\}.\] 
\item Consider the following special case of~\eqref{prob:ro_ddu} with \(\mb{x} \in \mathbb{R}^{n}\) \({\mb{y} \in \mathbb{R}^{m}\;,\;z \in \mathbb{R}}\)
  \begin{equation*}
  \label{RO-SAT}
  \tag{RO-SAT}
\min_{\mb{x,y},z \geq 0} \left\{ - z \mid  z - \bs{\xi}^\top \mb{y} \leq 0, \;\; \forall \bs{\xi} \in \mathcal{U}(\mb{x}),\;\;
     \mb{x,y} \leq \mb{e}, \; \;
    -\mb{y} \leq -\mb{e} \right\},     
  \end{equation*}
where $\mathcal{U}(\mb{x}) = \left\{(\xi_1,\dots,\xi_m) \mid \xi_i \geq x_{i_1} \;,\; \xi_i \geq x_{i_2} \;,\; \xi_i \geq 1 - x_{i_3} \;,\; \xi_i \leq 1 \right\}$.

Note that the 3-SAT problem is embedded in this set.
\item By Lemma~\ref{lem:npc} (provided after these steps), the optimal value of~\eqref{RO-SAT} is $-m$, if and only if the 3-SAT problem has a solution.
\item Problem~\eqref{RO-SAT} is a special case of~\eqref{prob:ro_ddu} with polyhedral set $\mU(\mb{x})$.
\item Since the 3-SAT problem is NP-complete~\citep{cook1971complexity}, problem~\eqref{prob:ro_ddu} is also NP-complete.
\end{enumerate}
\end{proof}
%%%%%%%%%%%%%%%%%%%%
\begin{lemma}
  \label{lem:npc}
  The 3-SAT problem has a feasible solution \(\mb{x}\), if and only if problem~\eqref{RO-SAT} has an optimal value of at most \(-m\).
\end{lemma}

\begin{proof}
\noindent\((\implies)\) Suppose the 3-SAT problem has a feasible solution \(\mb{x}\). This means, \(\mb{x}\) has to satisfy
\[x_{i_1} + x_{i_2} + (1 - x_{i_3}) \geq 1 \;\; \forall i = 1, \dots, m.\]
Since \(\mb{x} \in \{0,1\}^{n}\), for each \(i\) at least one of \(x_{i_1}\), \(x_{i_2}\), or \(1-x_{i_3}\) must be equal to 1.
Now, consider the uncertainty set 
$$\mathcal{U}(\mb{x}) = \left\{(\xi_1,\dots \xi_m) \mid \xi_i \geq x_{i_1} \;,\; \xi_i \geq x_{i_2} \;,\; \xi_i \geq 1 - x_{i_3} \;,\; \xi_i \leq 1 \;\;\forall i = 1, \dots, m  \right\}.$$
Since at least one of \(x_{i_1}\), \(x_{i_2}\), or \(1-x_{i_3}\) equals 1, \(\xi_i\) satsifies \(\xi_i \geq 1\). This implies that \(\xi_i = 1 \;\forall i\) is the only point in $\mathcal{U}(\mb{x})$. 
For this uncertainty set, the feasible solution is \(\mb{x, y} = \mb{1},\; z = m\). This leads to the optimal solution \(-z \leq -m\) or \(z \geq m\).

\((\impliedby)\) Suppose~\eqref{RO-SAT} has an optimal solution $(\mb{x}^*,\mb{y}^*)$ with the objective value of \({-z^{*} \leq -m}\). We first show that strict inequality is not possible.
Assume \(-z^{*} < -m\). 
The constraints in~\eqref{RO-SAT} imply \(z^* - \bs{\xi}^\top\mb{y}^* \leq 0\), i.e., \(\bs{\xi}^\top\mb{y}^* \geq z^{*} >m\) \(\forall \bs{\xi} \in \mathcal{U}(\mb{x}^*)\).
The constraints also imply \(y_i^* = 1 \; \forall i\). 
This means that \(\sum_{i=1}^{m} \xi_i > m\) \( \forall \bs{\xi} \in \mathcal{U}(\mb{x}^*)\).
However, the construction of the uncertainty set yields \(\xi_i \leq 1\). 
This leads to a contradiction, because \(\sum_{i=1}^{m} \xi_i \not> m\),  and hence \(-z^{*} = -m\).
Thus, \(\bs{\xi}^\top\mb{y}^* = m\) \(\forall \bs{\xi} \in \mathcal{U}(\mb{x}^*)\). 
Therefore, we can write
		\( \sum_{i=1}^{m} \xi_i = m  \;\; \forall \bs{\xi} \in \mathcal{U}(\mb{x}^*),\)
which implies
\( \min_{\bs{\xi} \in \mathcal{U}(\mb{x}^*)} \sum_{i=1}^{m}\xi_i	= m\).
However, since the uncertainty set implies \(\xi_i \leq 1 \;\forall i\), we can conclude that the sum can only be equal to $m$, if \(\xi_i = 1\;\; \forall i\).

We now show that this result implies for each \(i\) at least one of \(x_{i_1}^*\) or \(x_{i_2}^*\) or \((1-x_{i_3}^*)\) is equal to 1.
Suppose this is not true. This implies \(\exists i\) for which \(x_{i_1}^* < 1 \;,\; x_{i_2}^* < 1\) and \((1-x_{i_3}^*) <1 \). That means that we can construct \(\xi'_{i} = \max \{x_{i_1}^* , x_{i_2}^* , (1-x_{i_3}^*)\}\) which is an element of the uncertainty set and \(\xi'_{i} < 1\). 
However, this contradicts the result of \(\xi_i = 1 \;\; \forall i\).
Therefore, if \(z^* = m\), then we can find a feasible solution for the 3-SAT problem.
\end{proof}

Although problem~\eqref{prob:ro_ddu} is NP-complete, it can be reformulated as a bilinear or biconvex program, which may be solved by global optimization techniques~\citep[e.g.,][]{kolodziej2013global}.
For binary decision variables $\mb{x}$ influencing $\mathcal{U}(\mb{x})$, the problem~\eqref{prob:ro_ddu} can be reformulated as an MILP, using the \mbox{Big-M} method (see Section~\ref{sec:ext_gen}).
However, they suffer from weak relaxations.

%%%%%%%%%%%%%%%%%%%%%%%%%%%%%%%%%%%
\section{Structured Uncertainty Sets} 
\label{sec:red_ddu}
The weak numerical performance of \mbox{Big-M} linearization can be overcome, if the decision
$\mb{x}$ plays a decisive role in governing the elements of the uncertainty set.
Specifically, if the effect of $\mb{x}$ on the uncertainty set constraints can be modeled by penalizing the objective coefficients,
then the number of constraints in the robust counterpart can be reduced.
Here, we discuss the setting where $\mb{x}$ controls the upper bounds of the
uncertain variables.
This mechanism can be expressed in the set:
\begin{equation*}
\textbf{\(\overline{\bs{\Pi}}\)-Uncertainty:}  \quad \mU^{\overline{\bs{\Pi}}}(\mb{x}) = \left\{\bs{\xi} \mid \mb{D}\bs{\xi} \leq \mb{d},\;\; \bs{\xi} \leq \mb{v} + \mb{W}(\mb{e - x}),\;\;\bs{\xi} \geq \mb{0} \right\}.\label{prob:pibar} 
\end{equation*}
Here, $\mb{D}\in\mathbb{R}^{m\times n}$ is a coefficient matrix, \(\mb{d} \in \mathbb{R}^{m}\) is the RHS vector, $\mb{v} \in \mathbb{R}^n_{+}$ are the
minimum upper bounds, and $\mb{W} = diag(\mb{w}) \in \mathbb{R}^{n \times n}_{+}$ (a diagonal matrix) are the incremental
upper bounds, which apply when reduction is not applied. 
For \(\mU^{\overline{\bs{\Pi}}}\), the influence variable is \(\mb{x} \in \{0,1\}^n\).
The decision dependence in $\mU^{\overline{\bs{\Pi}}}$ affects the upper bounds on each
uncertain component $\xi_i$. 
This means, if the problem allows influencing uncertainties, this set can
model \emph{proactive} uncertainty reduction.
One possible example is disaster planning, where a decision to reduce the
fragility of certain roads yields an improved worst-case outcome.
Another example is measurement applications, where a decision for additional
expenditure leads to increased accuracy.
We employed such a set in Example 2 and discuss it further
in the numerical application.

We now discuss how this structure can be leveraged to reformulate the original problem~\eqref{prob:ro_ddu}.
Note that the objective function remains unaffected by the definition of the uncertainty set, as does the first term of the constraint.
Therefore, we focus our discussion on the parts of the constraint in problem~\eqref{prob:ro_ddu}, that are affected by uncertainty.

\subsection{$\overline{\bs{\Pi}}$-Uncertainty}
For succinctness, this section provides a reformulation of the following linear constraint
\begin{equation}
\label{eq:C}\tag{LC}
\mb{y}^{\top}\bs{\xi} \leq b\;\; \forall \bs{\xi} \in \mU^{\overline{\bs{\Pi}}}(\mb{x}).
\end{equation}
To satisfy this constraint for all \(\bs{\xi} \in \mU^{\overline{\bs{\Pi}}}(\mb{x})\), the uncertain LHS needs to be replaced by its maximum over the set. 
For this, consider the following two problems:

\begin{minipage}{0.5\textwidth}
\begin{equation}
  \label{prob:P}
  \tag{P}
  \begin{aligned}
    h(\mb{x,y}&) =\\
    \max_{\bs{\xi}}\;&\;\mb{y}^{\top}\bs{\xi}\\
    \text{s.t.}\;&\; \mb{D}\bs{\xi} \leq \mb{d} \\
    &\bs{\xi} \leq \mb{v} + \mb{W}(\mb{e - x}) & :\bs{\pi}(\mb{x,y})\\
    &\bs{\xi} \geq \mb{0},
  \end{aligned}
\end{equation}
\end{minipage}
\hspace{3mm}
\begin{minipage}{0.4\textwidth}
\begin{equation}
  \label{prob:Pp}
  \tag{P'}
\hspace{-5mm}
  \begin{aligned}
    \bar{h}(\mb{x,y}&) = \\
    \max_{\bs{\xi,\zeta}}\;&\;(\mb{y} - \overline{\bs{\Pi}}\mb{x})^{\top}\bs{\xi} + \mb{y}^{\top}\bs{\zeta}\\
    \text{s.t.}\;&\; \mb{D}\bs{\xi} +  \mb{D}\bs{\zeta} \leq \mb{d} \\
    &\bs{\xi} \leq \mb{W}\mb{e}\\
    &\bs{\zeta} \leq \mb{v}\\
    &\bs{\xi},\bs{\zeta} \geq \mb{0},
  \end{aligned}
\end{equation}
\end{minipage}
where in problem~\eqref{prob:P}, \(\bs{\pi}(\mb{x,y})\) denotes the corresponding dual variable. 
Problem~\eqref{prob:P} maximizes the LHS directly over $\mU^{\overline{\bs{\Pi}}}(\mb{x})$.
However, the standard reformulation of this problem leads to bilinear terms.
To avoid them, we can leverage the structure of the uncertainty set and
formulate problem~\eqref{prob:P} as problem~\eqref{prob:Pp}.
Such a problem pair was also suggested in the context of stochastic network interdiction~\cite{cormican1998stochastic}.
Proposition~\ref{prop:prop1} uses the duals of~\eqref{prob:P} and~\eqref{prob:Pp} to prove that they have the same objective value at optimality.
Formulating problem~\eqref{prob:Pp} requires the use of matrix \(\overline{\bs{\Pi}} = diag(\overline{\bs{\pi}})\).
Here, \(\overline{\bs{\pi}}\) is a component-wise upper bound of the optimal value of the dual variable \(\bs{\pi}(\mb{x,y})\) for all \(\mb{x,y}\). 
Note that the matrix \(\overline{\bs{\Pi}}\) is similar to \(M\) of the \mbox{Big-M} method in that it estimates an upper bound to the dual variables.
We provide a method to estimate \(\overline{\bs{\pi}}\) in Proposition~\ref{prop:est_pi}.
The dual problems of~\eqref{prob:P} and~\eqref{prob:Pp} are given by:

\begin{minipage}{0.48\textwidth}
  \begin{equation}
    \label{prob:D}
    \tag{D}
\hspace{-4mm}
    \begin{aligned}
      g(\mb{x,y}&) = \\
      \min_{\bs{\pi},\mb{q}}\;&\; \mb{q}^{\top}\mb{d} + \bs{\pi}^{\top}\mb{v} + \bs{\pi}^{\top}\mb{W}(\mb{e - x})\\
      \text{s.t.}\;&\;\bs{\pi}^{\top} + \mb{q}^{\top}\mb{D} \geq \mb{y}^{\top}\\
      &\; \bs{\pi},\mb{q} \geq \mb{0},
    \end{aligned}
  \end{equation}
\end{minipage}
\hspace{2mm}
\begin{minipage}{0.43\textwidth}
  \begin{equation}
    \label{prob:Dd}
    \tag{D'}
\hspace{-3mm}
    \begin{aligned}
      \bar{g}(\mb{x,y}&) = \\
      \min_{\mb{r,s,t}}\;&\; \mb{t}^{\top}\mb{d} + \mb{r}^{\top}\mb{We} + \mb{s}^{\top}\mb{v}\\
      \text{s.t.}\;&\; \mb{s}^{\top} + \mb{t}^{\top}\mb{D} \geq \mb{y}^{\top}\\
      &\; \mb{r}^{\top} + \mb{t}^{\top}\mb{D} \geq \mb{y}^{\top} - \mb{x}^{\top}\overline{\bs{\Pi}}\\
      &\; \mb{r,s,t} \geq \mb{0}.
    \end{aligned}
  \end{equation}
\end{minipage}

\begin{proposition}
  \label{prop:prop1}
Given a binary $\mb{x}$, if the set \(\mU^{\overline{\bs{\Pi}}}(\mb{x})\) is nonempty and $\mb{v,W}\ge0$, then for all \(\mb{y}\):
\[h(\mb{x,y}) = \bar{h}(\mb{x,y}).\]
\end{proposition}
\begin{proof}
Strong duality warrants the equalities \(g(\mb{x,y}) = h(\mb{x,y})\) and \(\bar{g}(\mb{x,y}) = \bar{h}(\mb{x,y})\).
In the following, we also refer to the optimal objective values of the dual problems as \(h(\mb{x,y})\) and \(\bar{h}(\mb{x,y})\).
Let \((\bs{\pi},\mb{q})\) be an optimal solution to~\eqref{prob:D}. 
Furthermore, let \({(\mb{r} = \bs{\pi} - \bs{\Pi}\mb{x}, \mb{s} = \bs{\pi},\mb{t} = \mb{q})}\) with \(\bs{\Pi} = diag(\bs{\pi})\) be a potential feasible solution to~\eqref{prob:Dd}.
For these solutions, it follows that 
\(\mb{s}^{\top} + \mb{t}^{\top}\mb{D} = \bs{\pi}^{\top} + \mb{q}^{\top}\mb{D} \geq \mb{y}^{\top} \),  and \({\mb{r}^{\top} +\mb{t}^{\top}\mb{D} = \bs{\pi}^{\top}- \mb{x}^\top\bs{\Pi} + \mb{q}^{\top}\mb{D} \geq \mb{y}^{\top} - \mb{x}^{\top}\bs{\Pi}
\geq \mb{y}^{\top} - \mb{x}^{\top}\overline{\bs{\Pi}}}\).
Since \({\bs{\pi},\mb{q} \geq \mb{0}}\), and \(\mb{x}\) is binary, we obtain \({\mb{r,s,t} \geq \mb{0}}\).
This means \((\mb{r,s,t})\) is a feasible solution to problem~\eqref{prob:Dd}.
This yields
\begin{eqnarray*}
\bar{h}(\mb{x,y}) &\leq& \mb{q}^{\top}\mb{d} + \bs{\pi}^{\top}\mb{v} + (\bs{\pi} -
\bs{\Pi}\mb{x})^{\top}\mb{We}\\ 
&=& h(\mb{x,y}).
\end{eqnarray*}
For the converse, let \({(\mb{r,s,t})}\) be an optimal solution to~\eqref{prob:Dd}.
Consider \({(\bs{\pi} = \mb{s},\mb{q} = \mb{t})}\) to be a solution to~\eqref{prob:D}.
The feasibility of \({(\mb{r,s,t})}\) leads 
\({\bs{\pi}^{\top} + \mb{q}^{\top}\mb{D} = \mb{s}^{\top} +
\mb{t}^{\top}\mb{D} \geq \mb{y}^{\top}}\), and \({\bs{\pi} = \mb{s} \geq \mb{0},\mb{q} = \mb{t} \geq \mb{0}}\). 
Hence, \((\bs{\pi},\mb{q})\) is a feasible solution to~\eqref{prob:D},
resulting in
\begin{eqnarray*}
h(\mb{x,y}) &\leq&  \mb{t}^{\top}\mb{d} + \mb{s}^{\top}\mb{v} +
  \mb{s}^{\top}\mb{W}(\mb{e - x})\\
&=& \bar{h}(\mb{x,y}) + (\mb{s}- \mb{r})^{\top}\mb{We} -
  \mb{s}^{\top}\mb{Wx}.
  \end{eqnarray*}
In order to prove \(h(\mb{x,y}) \leq \bar{h}(\mb{x,y})\), it is required to prove
\((\mb{s}- \mb{r})^{\top}\mb{We} - \mb{s}^{\top}\mb{Wx} \leq 0\).
This can be expressed as \(\sum_{i} w_i(s_i - r_i - s_i x_i)\le 0\).
For all $i$ with \(x_i = 1\), it holds that \({w_i(s_i - r_i  - s_i x_i) = - w_i r_i \leq 0}\).

Consider now the set of all \(i\) with \(x_i = 0\), denoted by \(X_0\).
Problem~\eqref{prob:Dd} can be rewritten as two nested minimization problems, where the outer problem is over \(\mb{t}\) and $r_j,s_j$ with $j\notin X_0$ and the inner problem over $r_i,s_i$ with $i \in X_0$:
\begin{equation*}
    \begin{aligned}
      \bar{h}(\mb{x,y}) = &\min_{\mb{t},r_j,s_j,j \notin X_0}\; \mb{t}^{\top}\mb{d}  + \sum_{j \notin X_0} r_j w_j  + \sum_{j \notin X_0} s_j v_j  + l(\mb{t})\\
      &\qquad\;\; {\left. 
      \begin{aligned} 
      \text{s.t.} \;&\; s_j + \mb{t}^{\top}\mb{D}_{\bigcdot,j} \geq y_j \\
      &\; r_j + \mb{t}^{\top}\mb{D}_{\bigcdot,j} \geq y_j - \overline{\pi}_j \\
      &\; r_j,s_j \geq \mb{0}
      \end{aligned} \quad\quad\right \rbrace}\; \forall j \notin X_0.
      \end{aligned}
  \end{equation*}
The inner minimization is captured by the function \(l(\mb{t})\), which is given by
\begin{equation*}
  \begin{aligned}
     l(\mb{t}) = &\min_{r_i,s_i,\;i \in X_0} \; \sum_{i \in X_0}r_i w_i  + \sum_{i \in X_0} s_i v_i\\
     &\qquad {\left.
     \begin{aligned}
      \text{s.t.}\;& s_i + \mb{t}^{\top}\mb{D}_{\bigcdot,i} \geq y_i \\ 
      &\; r_i + \mb{t}^{\top}\mb{D}_{\bigcdot,i} \geq y_i \\
      &\; r_i,s_i \geq \mb{0}
    \end{aligned} \quad\quad\right \rbrace} \forall \; i \in X_0.
    \end{aligned}
  \end{equation*}
Note that in this inner minimization problem, the same constraints act on \(s_i\) and \(r_i\). 
Since \(w_i\) and \(v_i\) are nonnegative, there exist optimal solutions \(s_i\) and \(r_i\) that are equal and set to their lower bounds \({s_i=r_i=\max\{y_i -  \mb{t}^{\top}\mb{D}_{\bigcdot,i}, 0\}}\). 
Therefore,  \({\sum_{i \in X_0} s_i w_i - r_i w_i  = 0}\), which means \({h(\mb{x,y}) \leq \bar{h}(\mb{x,y})}\). 
\end{proof}

Using Proposition~\ref{prop:prop1} and problem~\eqref{prob:Dd}, the constraint~\eqref{eq:C} can be reformulated as
\begin{equation*}
  \begin{aligned}
    &\; \mb{t}^{\top}\mb{d} + \mb{r}^{\top}\mb{We} + \mb{s}^{\top}\mb{v} \leq b\\
    &\; \mb{s}^{\top} + \mb{t}^{\top}\mb{D} \geq \mb{y}^{\top}\\
    &\; \mb{r}^{\top} + \mb{t}^{\top}\mb{D} \geq \mb{y}^{\top} - \mb{x}^{\top}\overline{\bs{\Pi}}\\
    &\; \mb{r,s,t} \geq \mb{0}.
  \end{aligned}
\end{equation*}
Note that this reformulation does not contain any bilinear terms and includes fewer constraints than the standard \mbox{Big-M} formulations.
Additionally, Proposition~\ref{prop:prop1} allows us to replace \(h(\mb{x,y})\) with \(\bar{h}(\mb{x,y})\). 
This is important because \(\bar{h}(\mb{x,y})\) is convex in \((\mb{x,y})\).
Therefore, cut generation algorithms can be used to solve this problem which is not possible for the original problem with the constraint~\eqref{eq:C}.
In the following, we discuss the matrix \(\overline{\bs{\Pi}}\).

\noindent\textbf{Estimation of \(\overline{\bs{\Pi}}\)} $\quad$
The following proposition sheds light on how to estimate \(\overline{\bs{\Pi}}\). 
\begin{proposition}
  \label{prop:est_pi}
If \(\mb{D}\) and $\mb{y}$ are element-wise nonnegative, then ${\pi_i(\mb{x,y}) \leq y_i}$ $\forall \mb{(x,y)}$ for constraint~\eqref{eq:C} under the uncertainty set $\mathcal{U}^{\overline{\bs{\Pi}}}$. 
\end{proposition}
\begin{proof}
  Consider the following problem for some index \(i\)
  \begin{equation}
    \label{P10s}
    \begin{aligned}
      F(\theta) = \max_{\bs{\xi}} \;&\;\mb{y}^{\top}\bs{\xi}\\
      \text{s.t.} \;&\;\mb{D}\bs{\xi} \leq \mb{d} &:\mb{q}\\
      &\; \bs{\xi} \leq \mb{v} + \mb{W}(\mb{e} - \mb{x}) + \theta \mb{e}_i &:\bs{\pi}\\
      &\; \bs{\xi} \geq \mb{0}.
    \end{aligned}
  \end{equation}
  Let \(\bs{\xi}_0\) be the optimal solution at \(\theta= 0\) and the corresponding optimal dual variables are \(\mb{q}_0\) and
  \(\bs{\pi}_0\).
  Let the optimal basis of the above problem be given by some matrix
  \(\mb{B}\).
  Since \(\bs{\xi}_0\) is the optimal solution, the vector of basic
  variables is given by \(\bs{\xi}_{0}^B = \mb{B}^{-1}\mb{b}\), where
  \(\mb{b}\) denotes the RHS vector of problem~\eqref{P10s}, i.e., \(\mb{b} =
  {[\mb{d}^{\top},\mb{v}^{\top} + {(\mb{e} - \mb{x})}^{\top}\mb{W}]}^{\top}\). 
  Assume that the solution is non-degenerate. This means \(\mb{B}^{-1}\mb{b} >
  0\). 
  Then for a small enough change in \(\mb{b}\), the optimal basis does not
  change. 
  If it is degenerate, then $\mb{b}$ can be perturbed by a small $\epsilon$ to obtain a non-degenerate solution, which only marginally changes the objective~\citep[see, e.g.,][]{bertsimas1997introduction}.

  When \(\theta >0\) is small enough, the basis matrix does not change.
  This means that both solutions (corresponding to \(\theta = 0\)
  and \(\theta > 0\)) have the same dual variables because the dual variables
  do not depend on the RHS vector.
  This means
  \begin{eqnarray*}
    F(\theta) - F(0) = &\;&\bs{\pi}_0^{\top}\mb{v} + \bs{\pi}_0^{\top}\mb{W}(\mb{e} - \mb{x}) + \theta \bs{\pi}_0^{\top}\mb{e}_i + \mb{q}_0^{\top} \mb{d} - \bs{\pi}_0^{\top}\mb{v} -\bs{\pi}_0^{\top}\mb{W}(\mb{e} - \mb{x}) - \mb{q}_0^{\top} \mb{d}\\
                         =&\;&\theta \bs{\pi}_0^{\top}\mb{e}_i,
  \end{eqnarray*}
  which represents the change in the objective value.
  Let \(\bs{\xi}_0\) be the optimal solution of the problem with \(\theta =
  0\) and \(\bs{\xi}_{\theta}\) be the  optimal solution of problem with \(\theta > 0\).
  Then the change in the objective value is 
  \[\theta \bs{\pi}_0^{\top}\mb{e}_i = \mb{y}^{\top} \bs{\xi}_{\theta} - \mb{y}^{\top} \bs{\xi}_0.\]

  \noindent Using Lemma~\ref{lem:ecbnd}, we can state that
  \begin{eqnarray*}
    \theta \bs{\pi}_0^{\top}\mb{e}_i &=& \mb{y}^{\top} \bs{\xi}_{\theta} - \mb{y}^{\top} \bs{\xi}_0\\
                                     &\leq&  \mb{y}^{\top} \bs{\xi}_{0} + \theta \mb{y}^{\top}\mb{e}_i - \mb{y}^{\top} \bs{\xi}_0 \\
                                     &=& \theta \mb{y}^{\top}\mb{e}_i.
  \end{eqnarray*}
  This implies that \({\pi}_{0,i} \leq {y}_i \;\;\forall i\).
\end{proof}
\begin{corollary}
Proposition~\ref{prop:est_pi} allows the estimation of $\overline{\Pi}$ by
  \begin{equation}
    \label{bound}
    \begin{aligned}
      \overline{\pi}_i = \max_{\mb{y}}\;&\; \mb{y}^{\top}\mb{e}_i\\
      \text{s.t.}\;&\; (\mb{x},\mb{y}) \in Y\\
      &\; x_{i} \in \{0,1\},
    \end{aligned}
  \end{equation}
where set $Y$ denotes the remaining constraints of the original full problem.
\end{corollary}

\begin{lemma}
  \label{lem:ecbnd}
If the matrix \(\mb{D}\) is element-wise greater than \(0\), then \(\bs{\xi}_{\theta} \leq \bs{\xi}_0 + \theta \mb{e}_i.\)
\end{lemma}
\begin{proof}
Suppose this is not true, i.e., there exists at least one index \(k\) such that \(\xi_{\theta,k} > \xi_{0,k} + \theta e_{i,k}\).
In addition, it holds that for $\theta\ge0$, \(\mb{y}^{\top}\bs{\xi}_{\theta} > \mb{y}^{\top}\bs{\xi}_0\).\\
If \(k \neq i\), then \(\bs{\xi}_{\theta} \leq \mb{v} + \mb{W}(\mb{e} - \mb{x})\), which suggests \(\bs{\xi}_{\theta}\) to be feasible for the problem with \(\theta = 0\).
This would contradict the optimality of \(\bs{\xi}_0\).\\
If \(k = i\), then \(\xi_{\theta,i} > \xi_{0,i} + \theta\).
However this results in
  \(\bs{\xi}_0 < \bs{\xi}_{\theta} - \theta \mb{e}_i \leq \mb{v} + \mb{W}(\mb{e} - \mb{x})\).
Since \(D (\bs{\xi}_{\theta} - \theta \mb{e}_i) = \mb{D}\bs{\xi}_{\theta} - \theta \mb{D}\mb{e}_i \leq \mb{d} - \theta \mb{D} \mb{e}_i \leq \mb{d}\),  \(\bs{\xi}_{\theta} - \theta \mb{e}_i\) is a feasible solution to the problem with \(\theta = 0\). 
However, this indicates that \(\mb{y}^{\top}(\bs{\xi}_{\theta} - \theta \mb{e}_i) > \mb{y}^{\top}\bs{\xi}_0\) which also contradicts the optimality of \(\bs{\xi}_0\).
Therefore, we can conclude that  \(\bs{\xi}_{\theta} \leq \bs{\xi}_{0} + \theta \mb{e}_i\).
\end{proof}
This proposition allows us to estimate \(\overline{{\pi}}_i\) by setting it equal to the maximum value that \(y_i\) can take in the overall problem. 
In some cases, such as shortest path or facility location problems, this is straightforwardly estimated from the underlying model. 
With this, all components of the \dd problem with the polyhedral uncertainty set $\mathcal{U}^{\overline{\bs{\Pi}}}$ can be computed efficiently for practical size problems.
We now extend Proposition~\ref{prop:prop1} to more general uncertainty sets.
% 

%%%%%%%%%%%%%%%%%%%%%%%%%%%%%%%%%%%%%%%%
\subsection{Extension to conic sets}
\label{sec:ext_conic}
%%%%%%%%%%%%%%%%%%%%%%%%%%%%%%%%%%%%%%%% 
Given a cone \(\mathcal{K}\), the \dd uncertainty set $\mathcal{U}^{\overline{\bs{\Pi}}}(\mb{x})$ can be extended  to
\begin{equation*}
  \mU^{\mathcal{K}}(\mb{x}) = \left\{\bs{\xi} \mid \mb{d} - \mb{D} \bs{\xi} \in \mathcal{K}, \; \bs{\xi} \leq \mb{v} + \mb{W}(\mb{e} - \mb{x}), \bs{\xi} \geq \mb{0}\right\}.  
\end{equation*}
 Here \(\mb{d}\) and \(\mb{D}\) are coefficients and $\mb{v}$ and \(\mb{W} = diag(\mb{w})\) denote upper bounds to the uncertain component \(\bs{\xi}\).
 The objective is to reformulate the constraint 
 \(\mb{y}^{\top}\bs{\xi} \leq b\), \(\forall \bs{\xi} \in \mU^{\mathcal{K}}(\mb{x})\).
In order to satisfy this constraint for all \(\bs{\xi} \in \mU^{\mathcal{K}}(\mb{x})\), its LHS can be expressed with the following two problems: 
 
 \begin{minipage}{0.5\textwidth}
 \begin{equation}
   \label{cP}
   \tag{KP}
   \hspace{-4mm}
   \begin{aligned}
     h(\mb{x,y}&)  = \\
     \max_{\bs{\xi}}\;&\; \mb{y}^{\top}\bs{\xi} \\
     \text{s.t.}\;&\; \mb{d} - \mb{D} \bs{\xi} \in \mathcal{K}\\
     &\; \bs{\xi} \leq \mb{v} + \mb{W}(\mb{e} - \mb{x}) \;\;:\;\bs{\pi}(\mb{x},\mb{y})\\
     &\; \bs{\xi} \geq \mb{0},
   \end{aligned}
 \end{equation}
 \end{minipage}%
 \hspace{5mm}
 \begin{minipage}{0.4\textwidth}
 \begin{equation}
   \label{cPp}
   \tag{KP'}
   \hspace{-4mm}
   \begin{aligned}
     \bar{h}(\mb{x,y}&)  =\\
      \max_{\bs{\xi,\zeta}}\;&\; (\mb{y} - \overline{\bs{\Pi}}\mb{x})^{\top}\bs{\xi} + \mb{y}^{\top}\bs{\zeta} \\
     \text{s.t.}\;&\; \mb{d} - \mb{D} \bs{\xi} \in \mathcal{K}\\
     &\; \bs{\xi} \leq \mb{We}\\
     &\; \bs{\zeta} \leq \mb{v}\\
     &\; \bs{\xi,\zeta} \geq \mb{0}.
   \end{aligned}
 \end{equation}
 \end{minipage}

\noindent Here, \(\bs{\pi}(\mb{x,y})\) denotes the dual variable for the corresponding constraint.
Let \(\overline{\bs{\Pi}}\) be an element-wise upper bound on the dual variables \({\bs{\pi}}(\mb{x,y})\).
The following proposition shows that the problems~\eqref{cP} and~\eqref{cPp} have the same optimal objective value. 
\begin{proposition}
  \label{prop:conicsets}
If \(\;\forall \mb{x} \in \{0,1\}^n\) there exists a point in the relative interior of $\;\mU^{\mathcal{K}}(\mb{x})$ (Slater point) and $\mb{v,W}\ge0$, then for all \(\mb{x, y}\):
  \[h(\mb{x,y}) = \bar{h}(\mb{x,y}).\]
\end{proposition}
The proof of this proposition proceeds similar to that of Proposition~\ref{prop:prop1}.
It uses strong duality which holds due to Slater's condition.
The proof proceeds parallel to the polyhedral uncertainty set.

Using Proposition~\ref{prop:conicsets} and the dual problem of~\eqref{cPp}, the constraint~\eqref{eq:C} can be reformulated as 
\begin{equation*}
  \begin{aligned}
    &\; \mb{t}^{\top}\mb{d} + \mb{r}^{\top}\mb{We} + \mb{s}^{\top}\mb{v} \leq b\\
    &\; \mb{s}^{\top} + \mb{t}^{\top}\mb{D} \geq \mb{y}^{\top}\\
    &\; \mb{r}^{\top} + \mb{t}^{\top}\mb{D} \geq \mb{y}^{\top} - \mb{x}^{\top}\overline{\bs{\Pi}}\\
    &\; \mb{t} \in \mathcal{K}^*,
    \; \mb{r,s} \geq \mb{0},
  \end{aligned}
\end{equation*}
with the dual cone $\mathcal{K}^*$.
Note that this reformulation has only linear terms and, as we will see in Section~\ref{sec:ext_gen}, fewer constraints than the standard \mbox{Big-M} formulation, hence it is more suitable for larger sized problems. 
The proof of this formulation proceeds parallel to that of Proposition~\ref{prop:prop1}.

In summary, these results allow the modeling of uncertainty sets with reducible upper bounds. 
Such bounds motivate the notion of \emph{proactive uncertainty control}. 
It mitigates conservatism and better actualizes the tradeoff between cost of control and disadvantage of uncertainty, both of which are instrumental parts of many real-world applications.
Until now, we discussed the special polyhedral set $\mU^{\overline{\bs{\Pi}}}$. 
The following section provides a reformulation of problem~\eqref{prob:ro_ddu} under general polyhedral uncertainty sets.

%%%%%%%%%%%%%%%%%%%%%%%%%%%%%%%%%%%%%%%%%%%%%%%%%%%%%%%%%%%%%%%%%%%%%%%%%%%%%%%

\section{Extensions to General Polyhedral Sets}
\label{sec:ext_gen}
The previous section leveraged the specific structure of the uncertainty set to obtain smaller reformulations. 
The \mbox{Big-M} reformulation, however, has the advantage of not requiring any special set structure. 
For completeness and a comparison of formulation sizes, the following proposition reformulates problem~\eqref{prob:ro_ddu} for the general polyhedral set $\mU^{P}(\mb{x})$. 

\begin{proposition}
\label{prop:gen_poly}
If the uncertainty set \(\mathcal{U}_{i}(\mathbf{x})\) is a polyhedron as in $\mU^{P}(\mb{x})$ with \(\bs{D}_{i} \in \mathbb{R}^{m_{i} \times p }\), \(\mb{d}_{i} \in \mathbb{R}^{m_{i}}\), and \(\bs{\Delta}_{i}
\in \mathbb{R}^{m_{i} \times n}\) and if  \(\mb{x}\) is binary, then the robust counterpart of problem \eqref{prob:ro_ddu} is
\begin{equation*}
  \label{prob:ro_ddu_milp}
  \begin{aligned}
    &\;\min_{\mb{x,y,w},\bs{\pi}} \;\;\mb{c}^{\top}\mb{x} + \mb{f}^{\top}\mb{y}\\
    &\;{\left.
    \begin{aligned}
    \text{s.t.} &\;\;\mb{a}_i^{\top}\mb{x} + \bs{\pi}_i^{\top}\mb{d}_i +  \sum_{j=1}^{m_i}\sum_{k=1}^{n}\Delta_{ijk}w_{ijk} \leq b_i \\
    &\; \bs{\pi}_i^{\top}\mb{D}_i = \mb{y}^{\top} 
    \end{aligned} \quad\right \rbrace} && \forall i\\
    &\quad {\left.
    \begin{aligned}
    &\;\;\; w_{ijk} \leq M x_k,\;\; w_{ijk} \leq \pi_{ij} \\
    &\;\;\; w_{ijk} \geq \pi_{ij} - M(1 - x_k) \\
    &\;\;\; \pi_{ij} \geq 0,\; w_{ijk} \geq 0 
    \end{aligned}\quad\quad\quad\quad\quad\quad\;\right \rbrace} && \forall i,j,k\\
    &\quad \;\;\;\;\; \mb{x} \in \{0,1\}^n, 
  \end{aligned}
\end{equation*}
where $M$ is a sufficiently large number.
\end{proposition}
\begin{proof}
We consider two cases, namely:
\emph{Case 1:}  There exists a feasible solution \((\mathbf{x}, \mathbf{y})\) to~\eqref{prob:ro_ddu}.
Therefore, \(\mathbf{x}\) and \(\mathbf{y}\) must satisfy all constraints
$\mb{a}_{i}^\top \mathbf{x} + \bs{\xi}_{i}^\top \mathbf{y} \leq b_{i}$ $\forall \bs{\xi}_{i} \in \mathcal{U}_{i}(\mathbf{x})$ for all $i$.
This is equivalent to 
\begin{align}\label{pr1} 
	\mb{a}_{i}^\top \mathbf{x}+\max_{\bs{\xi}_{i} \in \mathcal{U}_{i}(\mathbf{x})}\bs{\xi}_{i}^\top \mathbf{y}  \leq b_{i} \quad \forall i.
\end{align}
If this problem is feasible and has a finite optimal solution, then by strong duality, the corresponding dual problem has the same objective value.
Problem~\eqref{pr1} can now be expressed as
\begin{equation}		
\label{dual_constr}
{\left.
\begin{aligned}
&\> \mb{a}_{i}^\top \mathbf{x}+\bs{\pi}_{i}^\top (\mb{d}_{i} + \mb{\Delta}_{i}\mathbf{x}) \leq b_{i}	 &\qquad \\
&\> \bs{\pi}_{i}^\top \mb{D}_{i} =\mathbf{y}^\top  \qquad &\\
&\> \bs{\pi}_{i} \geq \mb{0}&\qquad \\
\end{aligned}\right \rbrace} \forall i,
\end{equation}
where $\bs{\pi}_i \in \mathbb{R}^{m_{i}}$ is the dual variable for constraints corresponding to the uncertainty set \(\mU_i(\mb{x})\). Here \(m_i\) refers to the number of constraints in the set \(\mU_i(\mb{x})\).
Since the primal problem is feasible and finitely valued, there exists a \(\bs{\pi}_{i}\), for which the constraints~\eqref{dual_constr} are satisfied.
Therefore, the original problem~\eqref{prob:ro_ddu} can be written as
\begin{equation*}
\begin{aligned} 
\>&\>\min_{\bs{\pi}_{i},\mb{x},\mb{y}}  \mb{c}^\top \mathbf{x} + \mb{f}^\top \mathbf{y}\\
&\quad{\left.
\begin{aligned}
\text{s.t.} \>&\>\mb{a}_{i}^\top \mathbf{x}+ \bs{\pi}_{i}^\top \mb{d}_{i} + \bs{\pi}_{i}^\top \mb{\Delta}_{i}\mathbf{x} \leq b_{i}	&\qquad\\
&\> \bs{\pi}_{i}^\top \mb{D}_{i} =\mathbf{y}^\top   &\qquad \\
&\> \bs{\pi}_{i} \geq \mb{0} &\qquad 
\end{aligned} \right \rbrace} \forall i. \end{aligned}
\end{equation*}
Note the bilinear term in the first constraint.
By expanding the variable space, the $i$th constraint can be rewritten as 
\begin{align*}
&\> \mb{a}_{i}^\top \mathbf{x}+ \sum_{j=1}^{m_{i}} \pi_{ij}d_{ij} + \sum_{j=1}^{m_{i}} \sum_{k=1}^{n}\Delta_{ijk} w_{ijk} \leq b_{i}, \;\; \text{with } w_{ijk} = \pi_{ij}x_{k}.
\end{align*}
In the bilinear term, \(w_{ijk} = \pi_{ij}x_{k}\),  \(x_{k}\) is binary, allowing to rewrite the term as
\begin{equation*}
	 w_{ijk} \leq \pi_{ij}, \;\;  0 \leq w_{ijk} \leq M x_{k},  \;\; w_{ijk} \geq \pi_{ij} - M(1- x_k),
\end{equation*}
where $M$ is a sufficiently large constant.		
Consequently, the problem~\eqref{prob:ro_ddu} can be reformulated as
\begin{equation}
\label{prob:ro_ddu_milp2}
\begin{aligned}
&\;\min_{\mb{x,y}} \;\;\mb{c}^{\top}\mb{x} + \mb{f}^{\top}\mb{y}\\
&\;{\left.
\begin{aligned}
\text{s.t.} &\;\;\mb{a}_i^{\top}\mb{x} + \bs{\pi}_i^{\top}\mb{d}_i +  \sum_{j=1}^{m_i}\sum_{k=1}^{n}\Delta_{ijk}w_{ijk} \leq b_i \\
&\; \bs{\pi}_i^{\top}\mb{D}_i = \mb{y}^{\top} 
\end{aligned} \quad\right \rbrace} && \forall i\\
&\quad {\left.
\begin{aligned}
&\;\;\; w_{ijk} \leq M x_k,\;\; w_{ijk} \leq \pi_{ij} \\
&\;\;\; w_{ijk} \geq \pi_{ij} - M(1 - x_k) \\
&\;\;\; \bs{\pi}_{i} \geq \mb{0},\; w_{ijk} \geq 0 
\end{aligned}\quad\quad\quad\quad\quad\quad\;\right \rbrace} && \forall i,j,k\\
&\quad \;\;\;\;\; \mb{x} \in \{0,1\}^n. 
\end{aligned}
\end{equation}
	
%%%%%%%%%%%%%%%%	
\emph{Case 2:}	 Problem~\eqref{prob:ro_ddu} is infeasible.
Then the reformulation in~\eqref{prob:ro_ddu_milp2} is infeasible.
To show this, consider the original problem~\eqref{prob:ro_ddu}.\\
Suppose this problem is infeasible under the assumptions of Proposition~\ref{prop:gen_poly}.
This means that \({\forall \mb{x}: \exists \bs{\xi} \in \mathcal{U}(\mb{x})}\) such that \(\mb{a}_{i}^\top \mathbf{x} + \bs{\xi}_{i}^\top \mathbf{y} > b_{i}\).
Consequently, the constraint \(\mb{a}_{i}^\top \mathbf{x}+\max_{\bs{\xi}_{i} \in \mathcal{U}_{i}(\mathbf{x})}\bs{\xi}_{i}^\top \mathbf{y}  > b_{i} \) holds for at least one \(i\).
Using the dual of the inner  problem, the constraints can be written \(\forall \bs{\pi}_i\) as
\begin{equation}		
\label{Stt:1}
\begin{aligned}
&\> \mb{a}_{i}^\top \mathbf{x}+\bs{\pi}_{i}^\top (\mb{d}_{i} + \mb{\Delta}_{i}\mathbf{x}) > b_{i}	 &\\
&\> \bs{\pi}_{i}^\top \mb{D}_{i} =\mb{y}^{\top}&\\
&\> \bs{\pi}_{i} \geq \mb{0}.&\\
\end{aligned}
\end{equation}

Now, assume that the reformulation in~\eqref{prob:ro_ddu_milp2} is feasible.	
Given its constraints, there exists a binary vector \(\mb{x}\) and a vector \(\mb{w}\) such that \(w_{ijk} = \pi_{ij}x_k\). 
However, this implies a variable \(\bs{\pi}_i = (\pi_{i1},\pi_{i2}, \dots, \pi_{ik}, \dots,\pi_{im_i})\) that satisfies \(\bs{\pi}_{i}^\top \mb{D}_{i} =\mathbf{y}^\top\), \(\bs{\pi}_i \geq \mb{0}\) and
\[\mb{a}_{i}^\top \mathbf{x}+\sum_{j=1}^{m_{i}} \pi_{ij}d_{ij} + \sum_{j=1}^{m_{i}} \sum_{k=1}^{n}\Delta_{ijk} \pi_{ij}x_{k} \leq b_{i}.\]
This contradicts the earlier assertion in~\eqref{Stt:1} that there exist no such \(\bs{\pi}_i\).
\end{proof}

This proposition allows us to reformulate the original \dd RO problem as a mixed
integer linear program which can be solved for many realistic size problems
using off-the-shelf algorithms.
Such mixed integer reformulations can also be provided for general convex uncertainty sets~\citep{ben2015deriving}, which includes conic and budgeted structures.
Their proofs (not shown) proceed parallel to that of Proposition~\ref{prop:gen_poly}.

Note that problem~\eqref{prob:ro_ddu} has \(n\) binary and \(p\) continuous variables, along with \(m\) constraints.
The i\(^\text{th}\) uncertain \(\bs{\xi}_i\) lies in an uncertainty set with \(m_i\) constraints. 
Table~\ref{tab:size} presents the size of the reformulation under two settings: 
\begin{enumerate*}[label=(\roman*)]
\item \(\mb{x}\) is binary as in Proposition~\ref{prop:gen_poly} and 
\item \(x_i\) can take \(s\) possible values.
\end{enumerate*}
For the sake of clarity, we assume that \(m_i = K \;\forall i\), where \(K\) is some constant. 
Table~\ref{tab:size} shows that for (ii), the size of the reformulation increases rapidly with growing \({s}\).
In certain cases, it is possible to improve the \mbox{Big-M} reformulation by imposing mild assumptions, as we will discuss next.
%%%%%%%%%%%%%%%%%%%%%%%%%%%%%
%%%%%%%%%%%%%%%%%%%%%%%%%%%%%
\begin{table}[h!]
\begin{tabular}{|l|l|l|l|l|} \hline
  Nature of \(\mb{x}\) & Binary var. & Continuous var.  & Affine constr. & Sign constr. \\ \hline
  Binary  & \(n\)  & \(p + mK + nK\)  & \(m + mp + 3nK\) & \(mK(n+1)\) \\ \hline
  Finite valued & \((n+1)s\)  & \(p + mK \)  & \(m + mp + 2n\)  & \(mK (ns+1)\)  \\ 
  &&\;\;\;\(+ nmK(s+1)\)&\(\;\;\;\;+ nmK(3s+1)\)&\\
  \hline
\end{tabular}
\caption{\label{tab:size} Size of \mbox{Big-M} formulation of~\eqref{prob:ro_ddu} for $\mathcal{U}_i(\mb{x})$ with respect to (i) \(\mb{x} \in \{0,1\}^n\) and (ii) \(\mb{x} \in \mathbb{R}^n\) with \(x_i\) taking \(s\) possible values:
$\dim(\mb{y}) = p$, $K$ constraints in \(\mU_i(\mb{x})\), and  $m$ constraints in the complete problem.}
\end{table}

\subsection{Modified \mbox{Big-M} Reformulation}
Consider the uncertainty set \(\mU^{\text{P}}(\mb{x})\) to be expressed as 
\begin{align*}
& \quad \mU^{\text{P}}(\mb{x}) = \left\{\bs{\xi} \mid \mb{D}_{i\bigcdot}^{\top}\bs{\xi} \leq d_i + \sum_{j=1}^{n}\Delta_{ij} x_j, \;\; \forall i=1,\dots,m \right\}. 
\end{align*}
To overcome the poor numerical performance of standard \mbox{Big-M} reformulation due to its weak relaxations, we impose the mild assumption that all elements of the coefficient matrix \(\bs{\Delta}\) are non-negative. 
Proposition~\ref{prop:rhsu} reformulates constraint~\eqref{eq:C} for \(\mU^\text{P}(\mb{x})\) under this assumption.  
\begin{proposition}
\label{prop:rhsu}
  If \(\Delta_{ij} \geq 0 \;\;\forall i,j\), then the constraint~\eqref{eq:C} with the uncertainty set $\mU^{\text{P}}(\mb{x})$ and a large constant $M$ can be reformulated as 
  \begin{equation*}
\begin{aligned}
  &\; \sum_{i=1}^{m} \pi_i d_i + \sum_{i=1}^{m} \sum_{j=1}^{n} t_{ij} \leq b\\
  &\; \sum_{i=1}^{m} \pi_i D_{ij} = y_j ,&& \forall j\\
  &\; {\left.
  \begin{aligned}
  &\;t_{ij} \geq \pi_i\Delta_{ij} - M(1 - x_j) \\  
  &\; \pi_i \geq 0, \;\; t_{ij} \geq 0 
  \end{aligned} \quad \quad \quad \right \rbrace } && \forall i,j.
\end{aligned}
\end{equation*}
\end{proposition}
\begin{proof}
  \noindent The LHS maximization problem for the constraint~\eqref{eq:C} can be written as 
  \begin{equation*}
    \begin{aligned}
      \max_{\bs{\xi}}\;&\; \mb{y}^{\top}\bs{\xi} \\
      \text{s.t.} \;&\; \mb{D}_{i\bigcdot}^{\top}\bs{\xi} \leq d_i + \sum_{j=1}^{n}\Delta_{ij} x_j\;\;\forall i.
    \end{aligned}
  \end{equation*}
Using the dual of this problem, the original constraint \(\mb{y}^{\top}\bs{\xi} \leq b\;\;\forall \bs{\xi} \in \mU^{\text{P}}(\mb{x})\) can be written as 
\begin{equation}
\label{conset1}
\begin{aligned}
  &\; \sum_{i=1}^{m} \pi_i (d_i + \sum_{j=1}^{n}\Delta_{ij} x_j) \leq b\\
  &\; \sum_{i=1}^{m} \pi_i D_{ij} = y_j \;\;\forall j\\
  &\; \bs{\pi} \geq \mb{0}.
\end{aligned}
\end{equation}
The constraints in~\eqref{conset1} can be rewritten by expanding the variable space as 
\begin{equation}
\label{conset2}
\begin{aligned}
  &\;  \sum_{i=1}^{m} \pi_i d_i + \sum_{i=1}^{m}\sum_{j=1}^{n} t_{ij} \leq b\\
  &\; \pi_i \Delta_{ij} x_j \leq t_{ij} \;&&\;\forall i, j\\
  &\; \sum_{i=1}^{m} \pi_i D_{ij} = y_j \;&&\;\forall j\\
  &\; \bs{\pi} \geq \mb{0}.
\end{aligned}
\end{equation}
If there is a variable \(\bs{\pi}\) feasible for the set of equations given by~\eqref{conset1}, then we can find a feasible variable for~\eqref{conset2} by \(t_{ij} = \pi_i\Delta_{ij} x_j\). 
On the other hand, if there exists a feasible solution to~\eqref{conset2}, then it is also feasible for~\eqref{conset1}.
If \(x_j = 0\), then \(t_{ij} \geq 0\) and if \(x_j = 1\), then \(t_{ij} \geq \pi_i\Delta_{ij}\). 
This can be expressed as the following set of constraints
\begin{equation*}
0 \leq t_j \geq \pi_i \Delta_{ij} - M(1 - x_j).
\end{equation*}
which completes the proof.
\end{proof}

The Proposition~\ref{prop:rhsu} leverages the fact that the variable \(t_{ij}\) remains at its lower bound, making the upper bounding constraints from the \mbox{Big-M} linearization redundant.
However, if \(t_{ij}\) can be negative, the two lower bounding constraints are not sufficient. 
%%
%%
%%%%%%%%%%%%%%%%%%%%%%%%%%%%%%%%%%%%%%%%%%%
%%%%%
\begin{table}[h!]
\centering
\begin{tabular}{|c|K{5.5cm}|c|}
\hline
 Formulations & Problem  &\setlength{\arrayrulewidth}{.05em}\begin{tabular}{c}Variables \\\hline Constraints \end{tabular}\\ \hline
 \(\overline{\bs{\Pi}}\) & 
\begin{minipage}[c]{1.0 \linewidth}
\begin{equation*}
  \begin{aligned}
  \\[-8mm]
    &\; \mb{t}^{\top}\mb{d} + \mb{r}^{\top}\mb{We} + \mb{s}^{\top}\mb{v} \leq b\\[-1mm]
    &\; \mb{s}^{\top} + \mb{t}^{\top}\mb{D} \geq \mb{y}^{\top}\\[-1mm]
    &\; \mb{r}^{\top} + \mb{t}^{\top}\mb{D} \geq \mb{y}^{\top} - \mb{x}^{\top}\overline{\bs{\Pi}}\\[-1mm]
    &\; \mb{r,s,t} \geq \mb{0}.
  \end{aligned}
\end{equation*}
\end{minipage}  
&\setlength{\arrayrulewidth}{.05em}\begin{tabular}{c}
C: \(m + 2n\)
\\[2mm]\hline\\[-2mm]
A: \(1 + 2n\) \\ S: \(m + 2n\)
\end{tabular} \\ \hline
%%%%
\mbox{Big-M}              & 
 \begin{minipage}[c]{1.0 \linewidth}
 \begin{equation*}
\begin{aligned}
\\[-8mm]
	&\; \mb{t}^{\top}\mb{d} + \mb{s}^{\top}\mb{v} + \bs{s}^{\top}\mb{We} - \sum_{i}r_i \leq b\\[-3mm]
    &\; \mb{s}^{\top} + \mb{t}^{\top}\mb{D} \geq \mb{y}^{\top}\\[-1mm]
    &\; w_{i}s_{i} - M (1- x_i) \leq r_i \leq M x_i\\[-1mm]
    &\; r_i \leq w_i s_i\\[-1mm]
    &\; \mb{r,s,t} \geq \mb{0}.
\end{aligned}
\end{equation*}
\end{minipage}
&\setlength{\arrayrulewidth}{.05em}\begin{tabular}{c}
C: \(m + 2n\)
\\[2mm]\hline\\[-2mm]
A: \(1 +4n\) \\ S: \(m + 2n\)
\end{tabular} \\ \hline
\begin{tabular}{c}Modified\\Big-M \end{tabular}
& \begin{minipage}[c]{1.0 \linewidth}
\begin{equation*}
\begin{aligned}
\\[-3mm]
	&\; \mb{t}^{\top}\mb{d} + \mb{s}^{\top}\mb{v} + \mb{r}^{\top}\mb{e} \leq b\\[-1mm]
    &\; \mb{s}^{\top} + \mb{t}^{\top}\mb{D} \geq \mb{y}^{\top}\\[-1mm]
    &\; r_{i} \geq w_{i}s_{i} - M x_i\\[-1mm]
    &\; \mb{r,s,t} \geq \mb{0}.
\end{aligned}
\end{equation*}
\end{minipage}  
&\setlength{\arrayrulewidth}{.05em}\begin{tabular}{c}
C: \(m + 2n\)
\\[2mm]\hline\\[-2mm]
A: \(1 + 2n\) \\ S: \(m + 2n\)
\end{tabular} \\ \hline
\end{tabular}
\vspace{3mm}
\caption{\label{tab:form}Comparison of~\eqref{eq:C} reformulations for the set $\mathcal{U}^{\overline{\bs{\Pi}}}(\mb{x})$ (C: Continuous, A: Affine, S: Sign).}
\end{table}
%%%%%%%%%%%%%%%%%%%%%%%%%%%%%%%%
%
%%%%%%%%%%%%%%%%%%%%%%%%%%%%%%%%
In some cases, it is possible to reformulate the problem even if the RHS coefficients are negative. 
Consider the shortest path example presented in the introduction, which has constraints of the form \(\xi_e \leq 1-\gamma_e x_e\). 
Here, the coefficient \({\Delta_e = -\gamma_e}\) is negative. 
However, we can rewrite the constraint as 
\({\xi_e \leq (1-\gamma_e) + \gamma_e (1 - x_e)}\) and apply the \mbox{Big-M} linearization on the variable \((1 - x_e)\) instead of on \(x_e\). 
This substitution allows the use of the modified \mbox{Big-M} reformulation in more general settings.
We report the numerical performance of this approach in comparison with the earlier reformulations in Section~\ref{sec:num_exp}.
For a comparison, we reformulate the constraint~\eqref{eq:C} over the uncertainty set \(\mU^{\overline{\bs{\Pi}}}(\mb{x})\)  using all three presented techniques, namely
\begin{enumerate*}[label=(\roman*)]
\item {$\overline{\bs{\Pi}}$},
\item \mbox{Big-M}, and
\item Modified \mbox{Big-M}.
\end{enumerate*}
Table~\ref{tab:form} presents this comparison along with the corresponding problem sizes.
The sign constraints correspond to \({(\bigcdot \geq 0)}\), which are presented separately since they can be solved more efficiently.
It displays that the primary difference between the \mbox{Big-M} and the other two reformulations is the larger number of affine (linear) constraints.
To gain intuition and provide computational comparison between the different formulations, 
we extend the introductory example of Section~\ref{sec:intro} to a more detailed numerical experiment.

%%%%%%%%%%%%%%%%%%%%%%%%%%%%%%%%%%%%%%%%%%%%%%%%%%%%%%%%%%%%%%%%%%%%%%%%%%%%%%%

\section{Numerical Experiments}
\label{sec:num_exp}

Shortest path problems on networks constitute a general class of models, describing the most efficient connection between a source and target. 
Deterministic shortest routing problems can be solved with polynomial time algorithms~\citep{dijkstra1959note}.
However, this does not hold for uncertain arc lengths.
Past research on robust shortest path problems focused on scenario-based~\citep{yu1998robust}, cardinality~\citep{bertsimas2003robust}, and interval uncertainty~\citep{averbakh2004interval,zielinski2004computational}.
Despite a large body of literature, to the best of our knowledge, there is no work in the context of uncertainties that depend on decisions.
To this end, our goals are:
\emph{
\begin{enumerate}
\item Comparing the numerical performance of different robust formulations,
\item Measuring the benefit of proactive reduction as a function of size, budget, or cost of reduction,
\item Measuring the number of arcs in the shortest path as a function of size, budget, or cost,
\item Evaluating the price of robustness and the benefit of interacting with uncertainties, and
\item Comparing the average and worst-case cost of decision dependence for \(RO\) and \(SO\).
\end{enumerate}
}
Here, we aim to model challenges that arise, e.g., in scenario planning of natural disasters.
When sections of a transportation network are damaged, the actual travel times along arcs become uncertain. 
To plan for such a scenario, a \dd RO solution can determine the arcs which should be strengthened (by reducing uncertainty) in order to improve the performance in an actual disaster.
This strengthening incurs a fee.
This means that it is possible to mitigate the impact of a disaster by managing the damage of a few particular arcs.
Similarly, for transportation problems (e.g., air, ground), travel time can be improved by acquiring additional  traffic or weather information on segments of the network.

To illustrate this setting, we discuss a problem on a graph \(G=(\mathcal{V},\mathcal{A},d(\bigcdot))\) for the set of nodes \(\mathcal{V}\), arcs \(\mathcal{A}\), and the distance function \(d(\bigcdot)\).
The objective is to find the shortest path from the source to the target node (\(s\ra t\)) when the actual realized distances from node \(i\) to \(j\) are uncertain and a function ${d}_{ij}(\bs{\xi}) = \left(1+ \frac{1}{2} \; \xi_{ij}\right)\bar{d}_{ij} $ of \(\bs{\xi}\).
The variable \(x_{ij}\) decides whether to reduce the maximum uncertainty in \({d}_{ij}\).
This inquiry comes at a cost \(c_{ij}\), which can be motivated as an investment in road improvement and is imposed on travelers via taxes or tolls.
The parameter $\bs{\xi}$ resides in a cardinality constrained uncertainty set with reducible upper bounds.
The complete problem is given by
\begin{equation}
  \label{prob:sp}
  \tag{SP}
  \begin{aligned}
    \min_{\mb{x},\mb{y}} \max_{\bs{\xi} \in \mU^{SP}(\mb{x})} \;&\;  \sum_{(i,j) \in \mathcal{A}} c_{ij} x_{ij} + \sum_{(i,j) \in \mathcal{A}} d_{ij}(\bs{\xi})y_{ij}\\
    \text{s.t.}\;&\; \mb{x} \in X \subseteq \{0,1\}^{|\mathcal{A}|},\;\mb{y} \in Y,
  \end{aligned}
\end{equation}
where \(y_{ij}\) decides whether the arc \((i,j)\) lies in the shortest path.
$X$ denotes any constraints on $\mb{x}$ and \(Y\) the set of routing constraints.
The uncertainty set is given by
\begin{equation*}
  \mU^{SP}(\mb{x}) = \left\{\bs{\xi} \mid 
    \sum_{(i,j) \in \mathcal{A}} \xi_{ij} \leq \Gamma,\; \xi_{ij} \leq\; 1 - \gamma_{ij} x_{ij} ,\;\xi_{ij} \geq 0\;\;\forall (i,j) \in \mathcal{A} \right\}.
\end{equation*}
We solve problem~\eqref{prob:sp} using the three different formulations: 
\begin{enumerate*}[label=(\roman*)]
\item  \(\overline{\bs{\Pi}}-\)formulation from Proposition~\ref{prop:prop1},
\item  standard  \mbox{Big-M} formulation, and
\item  Modified \mbox{Big-M} formulation from Proposition~\ref{prop:rhsu}.
\end{enumerate*}
%
%%%%%%%%%%%%%%%%%%%%%%%%%%%%%%%%%%%%
%%%%%%%%%%%
%
\begin{table}[h!]
\centering
\begin{tabular}{@{}|c|p{7.4cm}|p{25mm}|}
\hline
 Form. & \centering Problem  &\setlength{\arrayrulewidth}{.05em}\begin{tabular}{c}Variables \\\hline Constraints \end{tabular} \\ \hline
 \(\overline{\bs{\Pi}}\) & 
 \begin{minipage}[c]{1.0 \linewidth}
\begin{equation*}
\begin{aligned}
\\[-8mm]
\min_{\substack{\mb{x,y}\\ \mb{q,r},p}} \;& f(\mb{x,y}) \!+\! p \Gamma \!+\!\!\!\!\!\sum_{(i,j) \in \mathcal{A}}\!\!\! q_{ij}(1\! -\! \gamma_{ij}) \!+ \!\!\!\sum_{(i,j) \in \mathcal{A}}\!\! r_{ij} \gamma_{ij}\\
\text{s.t.} \;&\;  p + q_{ij} \geq \frac{y_{ij}d_{ij} - \overline{\pi}_{ij}d_{ij}x_{ij} }{2} \\[-3mm]
&\;  p + r_{ij} \geq \frac{y_{ij}d_{ij}}{2} \\[-1mm]
&\; p,q_{ij},r_{ij} \geq 0,\;
\mb{x},\mb{y} \in X \times Y. 
\end{aligned}
\end{equation*}
\end{minipage}
&\setlength{\arrayrulewidth}{.05em}\begin{tabular}{c}
B: \(2|\mathcal{A}|\) \\ C: \(2|\mathcal{A}|\!+\! 1\)
\\[2mm]\hline\\[-2mm]
A:\! \(|\mathcal{V}| \!+\! 2 |\mathcal{A}|\) \\ S:\! \(2|\mathcal{A}| + 1\)
\end{tabular} \\ \hline
 \mbox{Big-M}              & 
 \begin{minipage}[c]{1.0 \linewidth}
\begin{equation*}
\begin{aligned}
\\[-8mm]
\min_{\substack{\mb{x,y}\\ \mb{q,r},p}}&\> f(\mb{x,y}) + p \Gamma + \!\!\!\sum_{(i,j) \in \mathcal{A}} q_{ij} - \!\!\!\sum_{(i,j) \in \mathcal{A}}\gamma_{ij}r_{ij}\\[-2mm]
\text{s.t.} \>&\>  p + q_{ij} \geq \frac{d_{ij}y_{ij}}{2}\\[-1mm]
&\; 0 \leq r_{ij} \leq M x_{ij}\\[-1mm]
&\; q_{ij} - M(1 - x_{ij}) \leq r_{ij} \leq q_{ij}\\[-1mm]
&\> p, q_{ij},r_{ij} \geq 0,\; 
 \mb{x},\mb{y} \in X \times Y.
\end{aligned}
\end{equation*}
\end{minipage} 
&\setlength{\arrayrulewidth}{.05em}\begin{tabular}{c}
B: \(2|\mathcal{A}|\) \\ C: \(2|\mathcal{A}|\!+\! 1\)
\\[2mm]\hline\\[-2mm]
A: \(|\mathcal{V}| + 4 |\mathcal{A}|\) \\ S: \(2|\mathcal{A}| + 1\)
\end{tabular} \\ \hline
\begin{tabular}{c}Modified\\Big-M \end{tabular}
& \begin{minipage}[c]{1.0 \linewidth}
\begin{equation*}
\begin{aligned}
\\[-8mm]
\min_{\substack{\mb{x,y}\\ \mb{q,r},p}}\;&\;\; f(\mb{x,y}) + p \Gamma + \!\!\!\!\!\sum_{(i,j) \in \mathcal{A}}\!\! r_{ij} \!+ \!\!\!\!\! \sum_{(i,j) \in \mathcal{A}}\!\! q_{ij}(1 - \gamma_{ij})\\[-2mm]
\text{s.t.} \;&\; p + q_{ij} \geq \frac{d_{ij}y_{ij}}{2} \\[-1mm]
&\; r_{ij} \geq \gamma_{ij} - M x_{ij} \\[-1mm]
&\; p,  q_{ij},r_{ij} \geq 0,\;
\mb{x},\mb{y} \in X \times Y.
\end{aligned}
\end{equation*}
\end{minipage}  
&\setlength{\arrayrulewidth}{.05em}\begin{tabular}{c} 
B: \(2|\mathcal{A}|\) \\ C: \(2|\mathcal{A}|\!+\! 1\)
\\[2mm]\hline\\[-2mm]
A:\! \(|\mathcal{V}| + 2 |\mathcal{A}|\) \\ S:\! \(2|\mathcal{A}| \!+\! 1\)
\end{tabular} \\ \hline
\end{tabular}
\caption{\label{tab:formulations}Shortest path formulations for the set $\mU^{SP}(\mb{x})$ (B: Binary, C: Continuous, A: Affine, S: Sign).}
\end{table}
%%%%%%%%%%%%%%%%%
%%%%%%%%%%%%%%%%%
In Table~\ref{tab:formulations}, \(X\times Y\) denote the collection of both the shortest path and decision constraints.
Furthermore, \({f(\mb{x,y}) = \sum_{(i,j) \in \mathcal{A}} c_{ij} x_{ij} + \sum_{(i,j) \in \mathcal{A}} \bar{d}_{ij}y_{ij}}\) denotes the total cost of reduction and nominal length. 
Table~\ref{tab:formulations} shows that the difference between the \mbox{Big-M} formulation and the other two formulations lies in the number of affine (linear) constraints, as in Table~\ref{tab:form}. 
We now discuss the numerical experiments.

\emph{Experiment 1: Performance Comparison} \quad 
The numerical setup is as follows.
We randomly generate points on a \(100 \times 100\) area and connect them to
create a complete graph. 
The two furthest nodes constitute the source and destination. 
The final graph is selected after removing \(60\%\) of the longest arcs in order to avoid direct connections between the source and destination.
The uncertainty budget \(\Gamma\) is set to 2. 
The cost of reduction \(c_{ij} = c\) and the fraction of uncertainty reduced \(\gamma_{ij} = \gamma\) are \(1.0\) and \(0.2\), respectively.
For each size \(|\mathcal{V}| = \{50,75,\dots,300\}\), \(100\) random graphs are generated.
These values serve as an illustration of the qualitative comparison of the formulations.
In practical applications, they need to be estimated from the economical value of travel time ($d_{ij}$) relative to the per-trip tax burden for road investments ($c_{ij}$).
%%%%%%%%%%%%%%%%%%%
%
\begin{figure}[h!]
  \centering
  \includegraphics[width=80mm]{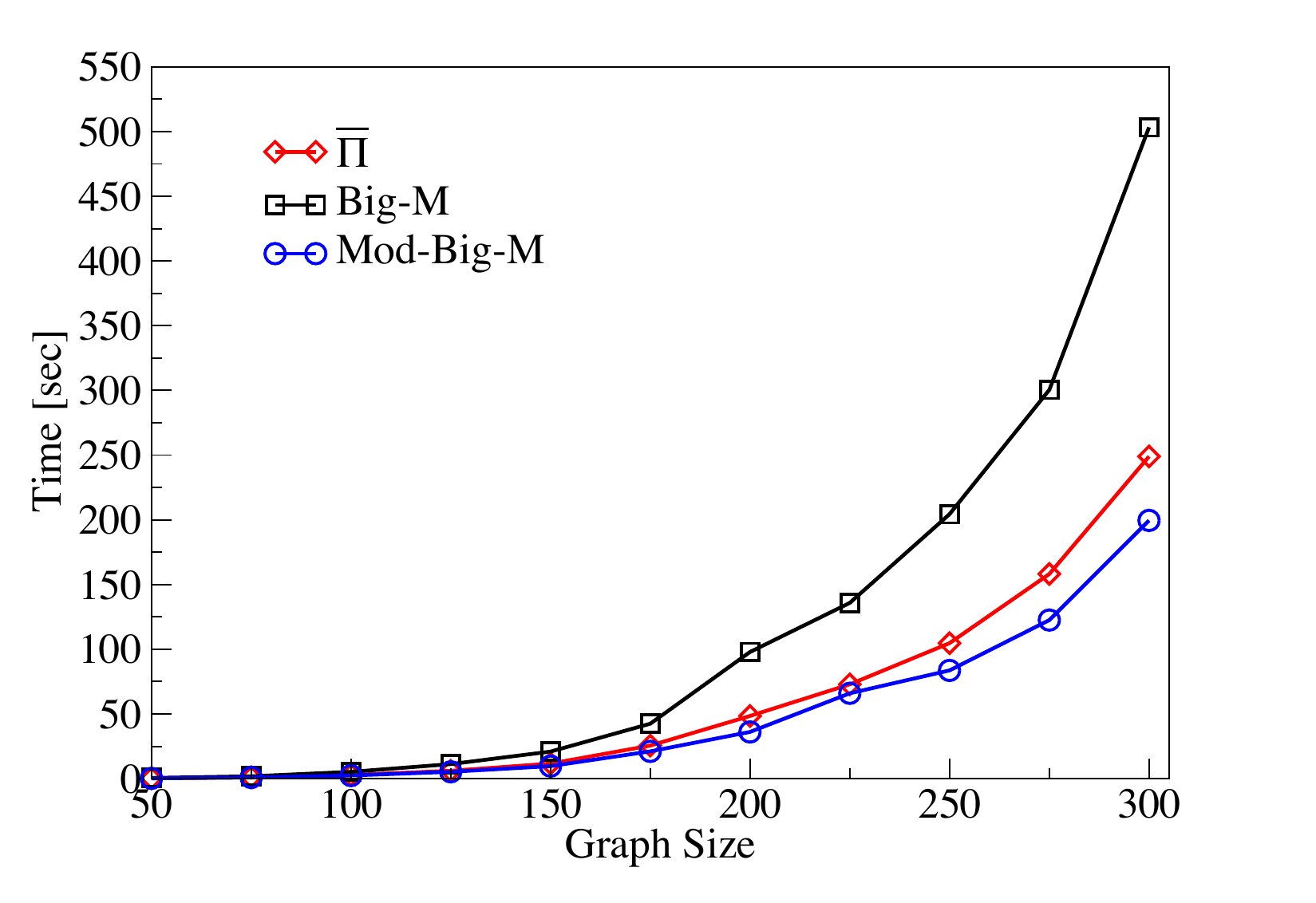}
  \vspace{-3mm}
  \caption{Comparison of median solution times of reformulations from Propositions~\ref{prop:prop1},~\ref{prop:rhsu}, and the standard \mbox{Big-M}.}
  \label{fig:solntime}
\end{figure}

To solve these problems, we used the Gurobi 7.0 solver on a commercially available computing unit with Intel Core i7 at 3.6~GHz.
The median computation times for different approaches and varying sizes are reported in Figure~\ref{fig:solntime}.
Note that all three methods lead to the same solution.
The observations from Figure~\ref{fig:solntime} can be summarized as follows.
\begin{itemize}
\item The time increases with growing \(|\mathcal{V}|\) for all formulations. 
However, the increase is less steep for the \(\overline{\bs{\Pi}}\) and the Modified \mbox{Big-M} formulation than for the \mbox{Big-M} formulation. 
\item The difference between the \mbox{Big-M} and the proposed formulations increases with growing \(|\mathcal{V}|\). 
  This highlights the advantage of the \(\overline{\bs{\Pi}}\) and Modified \mbox{Big-M} formulation for larger  graphs.
\item The median time of the Modified \mbox{Big-M} formulation is less than that of the \(\overline{\bs{\Pi}}\)-formulation. 
\end{itemize}
Figure~\ref{fig:solntime} highlights the benefits of using the proposed
formulations to solve such \dd optimization problems.
While the performance of the Modified \mbox{Big-M} and \(\overline{\bs{\Pi}}\) formulations are comparable over a broad range of network sizes, the subproblem in the \(\overline{\bs{\Pi}}\) reformulation is convex, which can be exploited by cut-generating methods, which may be computationally advantageous.
We also solved the \(\overline{\bs{\Pi}}\) formulation using a cut generation approach (not shown).
However, for this application, it converged slowly and required a sizable number of cuts. 

We now focus on analyzing how the solution changes as the parameters of the uncertainty set are varied. 
For this purpose, we introduce additional notation for observable quantities.

\noindent\emph{Notation for Observables. $\quad$} 
The number of arcs in the shortest path is $n^*$, which is a function of the budget $\Gamma$ and the level of uncertainty reduction $\gamma$.
These parameters create three scenarios:

\begin{enumerate}[label=(\emph{\roman*})]
\item  \emph{nominal} case, where no uncertainty is present, 
$n^*(\Gamma=0,\gamma=0)$;
\item  \emph{standard robust} case with no decision dependence, 
$n^*(\Gamma>0,\gamma=0)$; and
\item  \emph{\dd robust} case with uncertainty reduction 
${n^*(\Gamma>0,\gamma>0)}$, in which case $\tilde{n}$ is the number of arcs whose uncertainty was reduced.
\end{enumerate}
We also follow this notation for the optimal objective value $z^*$.
Consequently, the difference $\big(z^*(\Gamma>0,\gamma=0) - z^*(\Gamma=0,\gamma=0)\big)$ constitutes the \emph{price of robustness}, whereas the difference $\big(z^*(\Gamma>0,\gamma=0) - z^*(\Gamma>0,\gamma>0)\big)$ constitutes the \emph{benefit of interaction.}

There are four parameters that govern the effect of interactions with uncertainty: \(\gamma, |\mathcal{V}|, c,\) and  \(\Gamma\).
To evaluate their role and to infer the underlying mechanism, we devise four experiments by tuning across their range. 
Specifically, by adjusting one parameter while keeping the other three fixed, we explore four orthogonal settings.

In these experiments, the problem \eqref{prob:sp} is implemented on randomly generated graphs of $[20-50]$ nodes. 
This size is comparable to moderately sized transportation networks~ \citep{montemanni2005robust}. 
For each size, 2000 graphs are generated in a manner similar to the previous experiment.
We maintain these parameter values throughout the following experiments, except in those where their change is probed.
In the following, we discuss the four experiments.

\begin{figure}[h!]
	\centering	
	    \begin{minipage}[t]{75mm}
		\includegraphics[width=75mm]{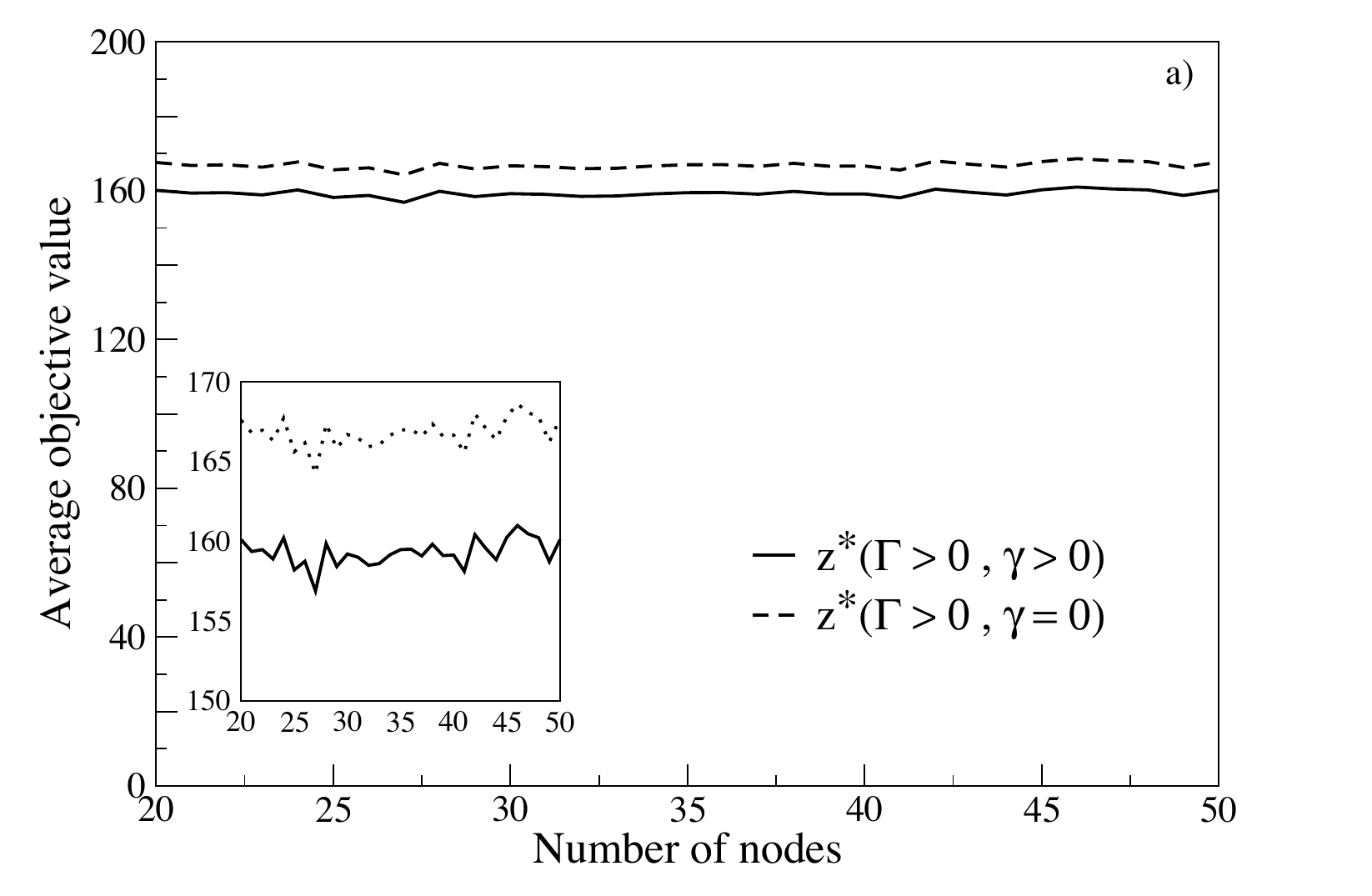}
	   \end{minipage}
	\hspace{-5mm}
  	   \begin{minipage}[t]{75mm}
		\includegraphics[width=75mm]{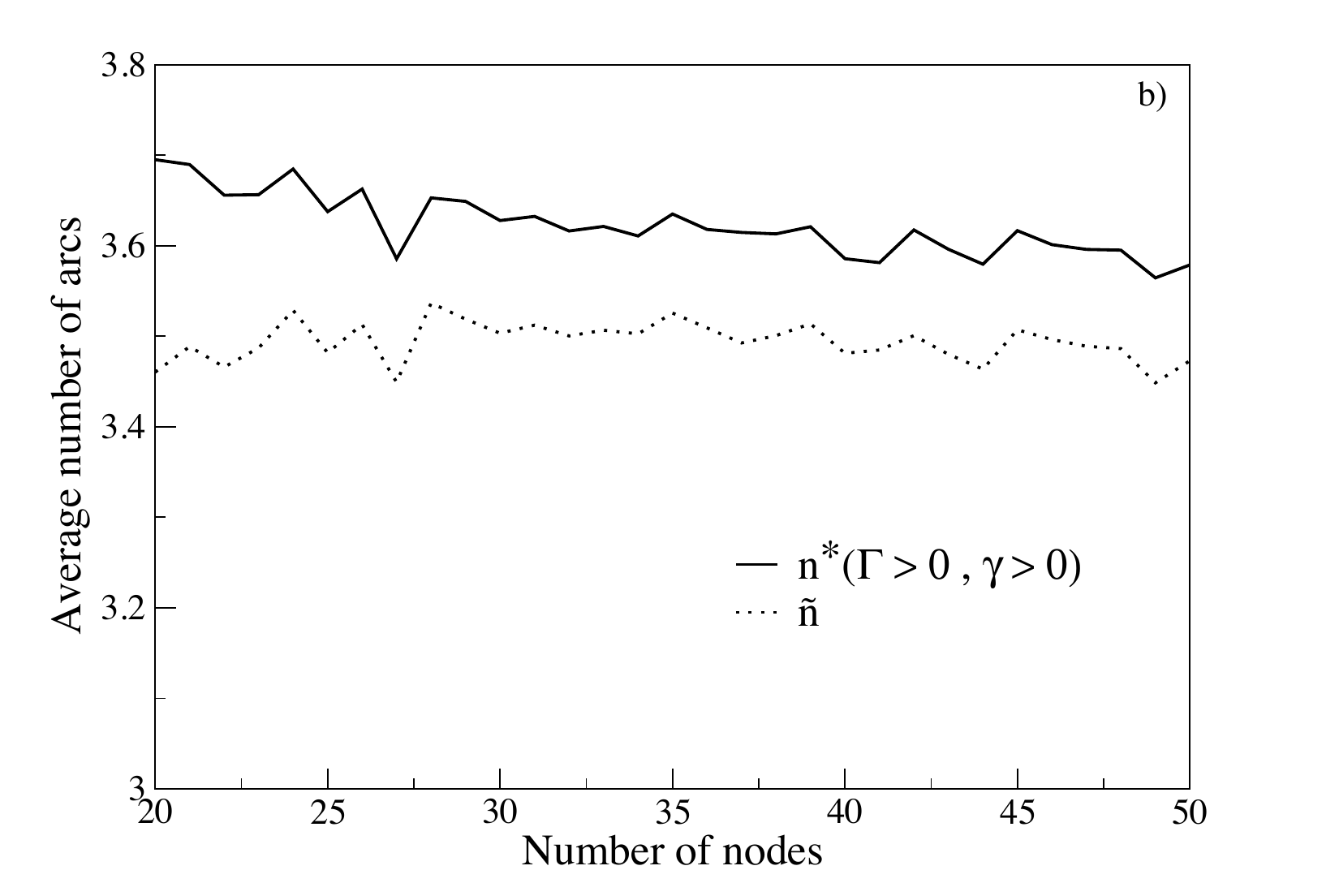}
	   \end{minipage}
\caption{Dependence on graph size \(|\mathcal{V}|\) for: a) average objective function and b) average number of arcs.  The inset is a magnification.}
\label{figr:12}	   
\end{figure}

%%%%%%%%%%%%%%%%%%%%%%%%%%%%%%%%%%%
\noindent\emph{Experiment 2: Uncertainty Reduction. $\quad$}
We compare $z^*$, when reduction is permitted ($\gamma>0$) or not ($\gamma=0$).
Figure~\ref{figr:12}a shows that $\gamma>0$ reduces $z^*$ (shorter paths),
which is independent of \(|\mathcal{V}|\).
The inset of Figure~\ref{figr:12}a is a magnification, displaying marginal fluctuations that stem from the random nature of graphs.

%%%%%%%%%%%%%%%%%%%%%%%%%%%%%%%%
\noindent\emph{Experiment 3: Graph Size. $\quad$}
We observe that not all arcs in the shortest path experience uncertainty reduction ($\tilde{n} <  n^*(\Gamma>0,\gamma>0)$), independent of \(|\mathcal{V}|\).
This is attributed to the non-zero \(c\).
We also observe that $z^*$ is independent of \(|\mathcal{V}|\), which can be explained by the fact that
\(|\mathcal{V}|\) only increases from \(20-50\) and ${n^*(\Gamma>0,\gamma>0)}$ does not change sizably over this range as such the effect on $z^*$ is undetectable. 
We expect \(n^*\) and $z^*$ to increase measurably when \(|\mathcal{V}|\) varies by a few orders of magnitude.
Larger experiments come at a significant computational burden and are outside the scope of this study.

Figure~\ref{figr:12}b illustrates the average $n^*(\Gamma>0,\gamma>0)$ and the average \(\tilde{n}\) for varying \(|\mathcal{V}|\).
We also observe a slight downward trend of $n^*(\Gamma>0,\gamma>0)$ with increasing \(|\mathcal{V}|\).
This is because the connectivity within a graph increases with \(|\mathcal{V}|\) as 
the number of arcs grows faster than the number of nodes, because in the experimental setup, only a fixed fraction of arcs are removed.
\begin{figure}[h]
	\centering
	    \begin{minipage}[t]{75mm}
		\includegraphics[width=75mm]{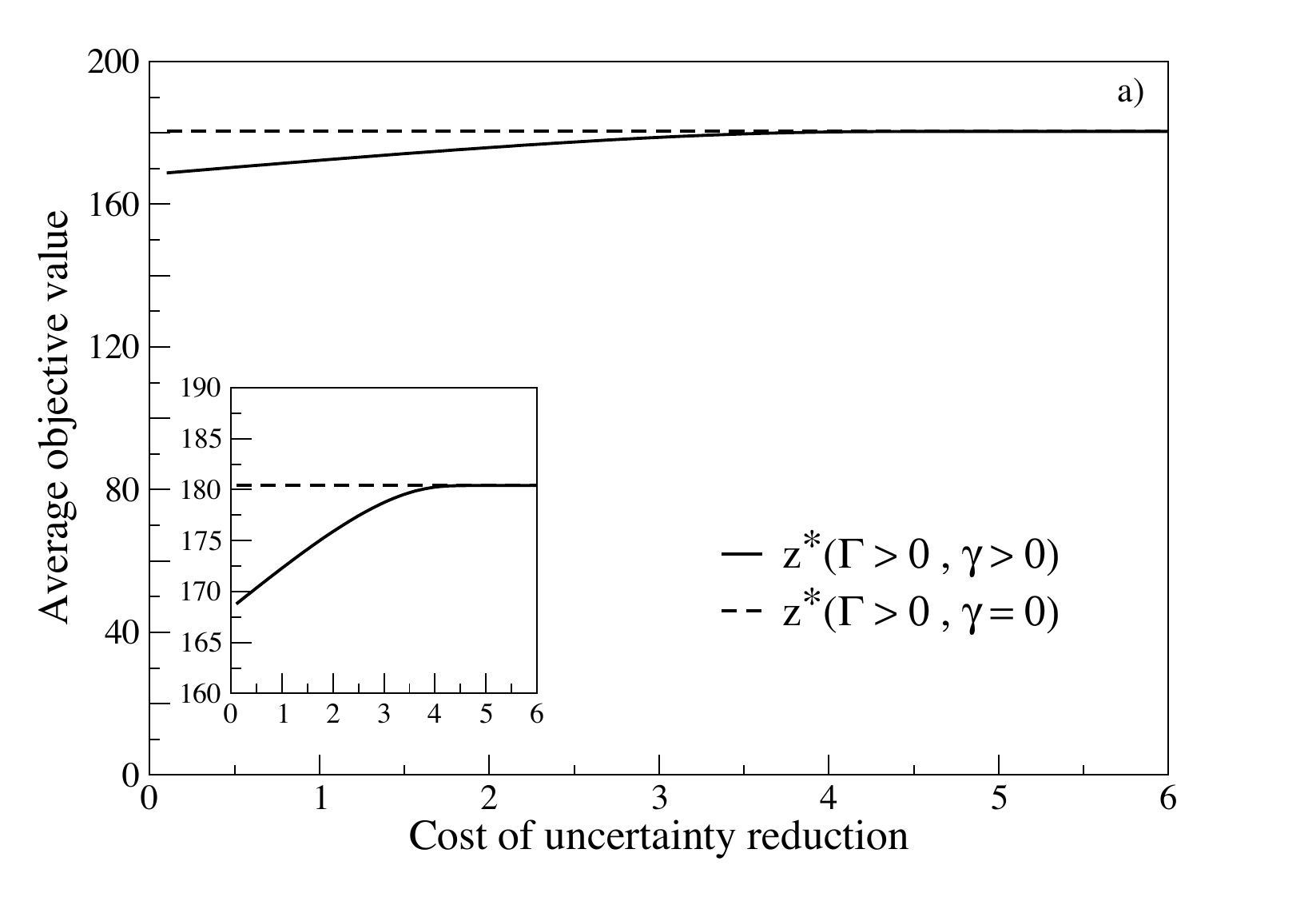}
	   \end{minipage}
	\hspace{-5mm}
  	   \begin{minipage}[t]{75mm}		
		\includegraphics[width=75mm]{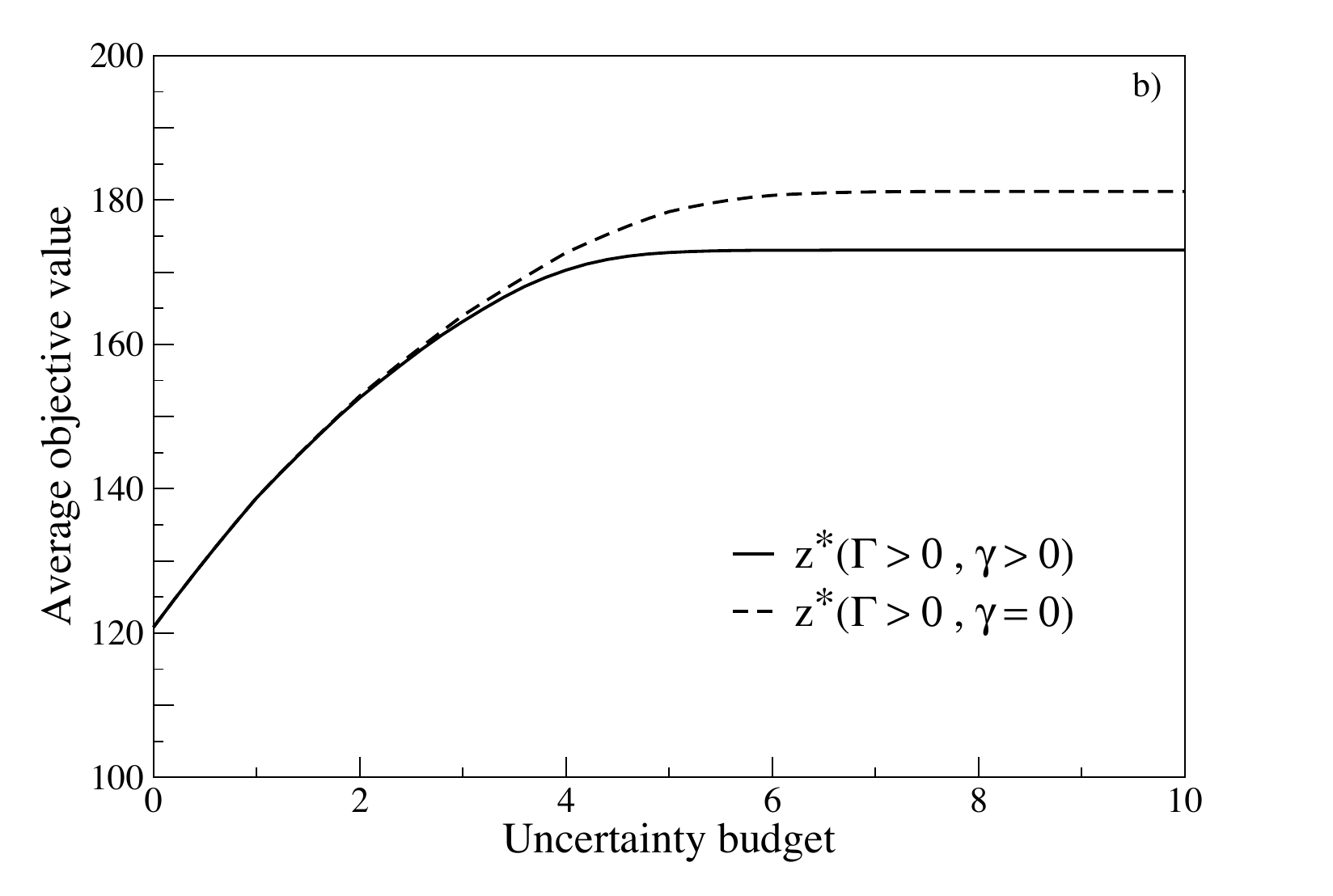}
	   \end{minipage}
\caption{Average objective value as a function of: a) cost of uncertainty reduction $c$ and b) maximum uncertainty $\Gamma$. The graph consists of \(|\mathcal{V}|=30\) nodes.}
\label{figr:35}
\end{figure}

\begin{figure}[h]
	\centering
	\begin{minipage}[t]{75mm}
		\includegraphics[width=75mm]{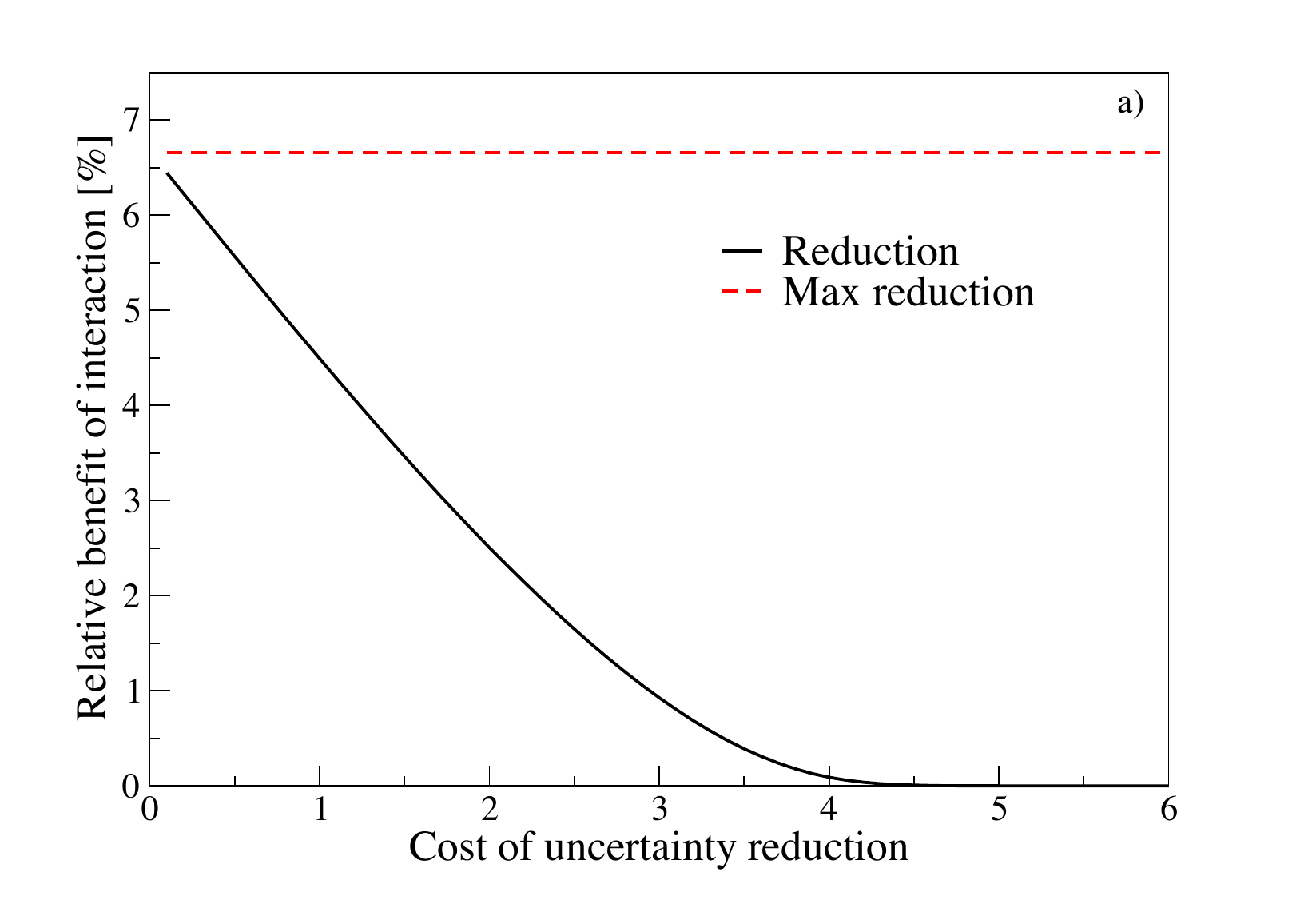}
	\end{minipage}
	\hspace{-5mm}
	\begin{minipage}[t]{75mm}		
		\includegraphics[width=75mm]{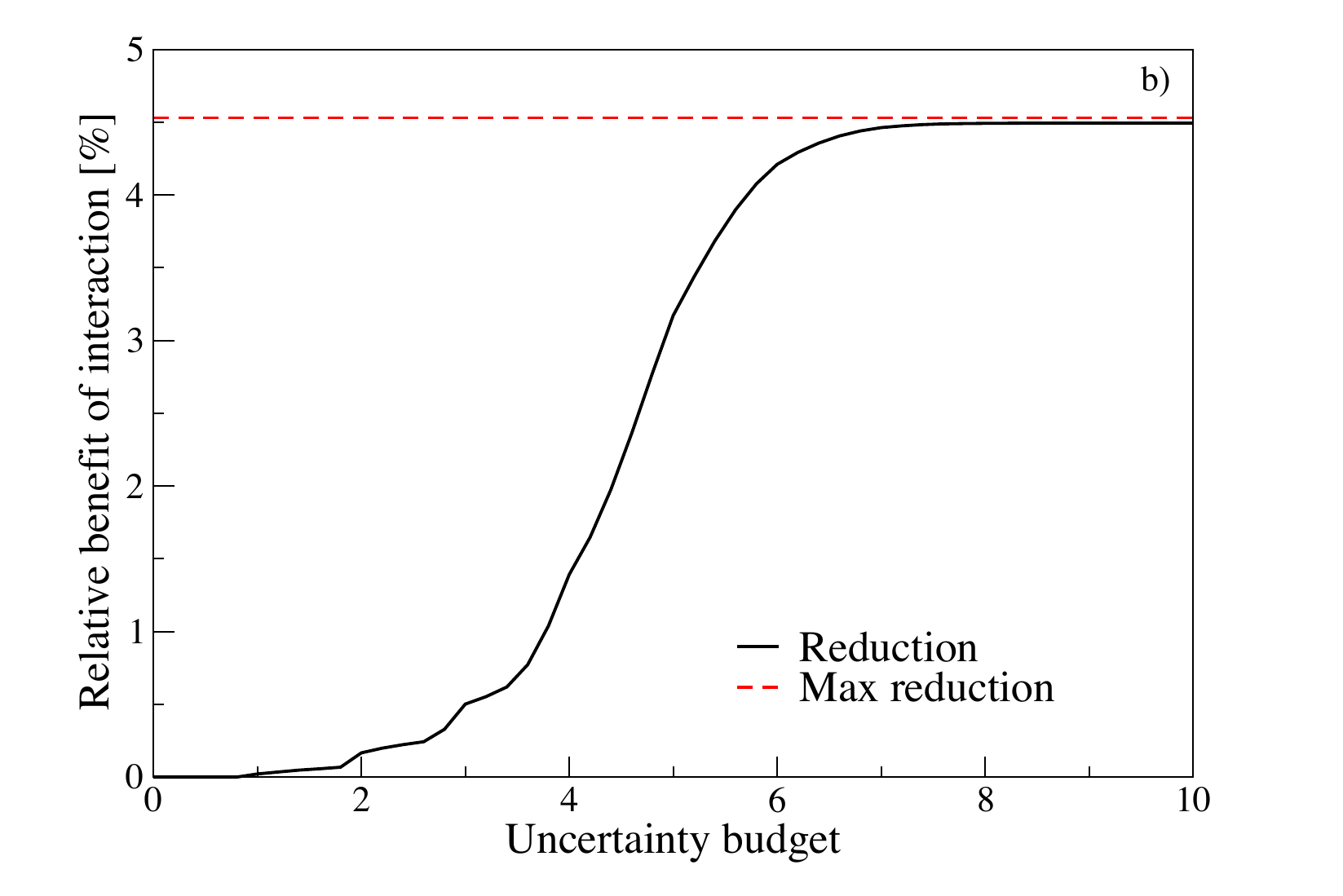}
	\end{minipage}
	\caption{Average relative benefit of interaction as a function of: a) cost of uncertainty reduction $c$ and b) maximum uncertainty $\Gamma$. The graph consists of \(|\mathcal{V}|=30\) nodes.}
	\label{figr:bor}
\end{figure}

\noindent\emph{Experiment 4: Cost of Uncertainty Reduction. $\quad$}
The reduction cost $c$ determines the trade-off between accepting the uncertainty level and its reduction.
It can be expected that an increasing $c$ marginalizes the benefits of reducing uncertainty.
This means that for a sufficiently low $c$, uncertainty can be reduced in every arc in the shortest path.
On the other hand, for high $c$, the opposite is true.
Figure~\ref{figr:35}a (\(|\mathcal{V}| = 30\; \text{and}\; \Gamma = 12\)) shows that for $c\le4$, the average \(z\) can be decreased.
However for large $c$, the high cost of reduction makes it disadvantageous to reduce uncertainty. 
%Benefit of reduction
The price of robustness (difference between the dotted line in Figure~\ref{figr:35}a and \(z^*(\Gamma = 0,\gamma = 0)\) in Figure~\ref{figr:35}b) is constant w.r.t. \(\gamma\) but changes with \(\Gamma\).
On the other hand, the benefit of interaction decreases with increase in \(c\), as can be observed in Figure~\ref{figr:bor}a.
Note that the maximum benefit of interaction is calculated by assuming uncertainty is reduced on all the arcs in the shortest path, at zero cost (\(c = 0\)).

%%%%%%%%%%%%%%%%%%%%%%%%%%%%%
\noindent\emph{Experiment 5: Uncertainty Budget. $\quad$}
\(\Gamma\) governs the number of arcs that can be affected by uncertainty.
Figure~\ref{figr:35}b shows that \(z^*\) increases gradually with $\Gamma$ until it reaches the level of the corresponding shortest path length affected by the relative uncertainty $(1+\frac{1}{2})$ and plateaus thereafter.
This is because increasing \(\Gamma\) beyond a certain point does not have any effect on \(n^*\), since all the arcs in the path are already uncertain and additional budget remains untapped.
Consequently, the price of robustness increases with \(\Gamma\) and plateaus beyond a certain \(\Gamma\) (not shown). 
An analogous behavior can be observed for the benefit of interaction, as shown in Figure~\ref{figr:bor}b. 
The maximum benefit is achieved at \(c = 0\).

\begin{figure}
	\centering
	    \begin{minipage}[t]{75mm}
		\includegraphics[width=75mm]{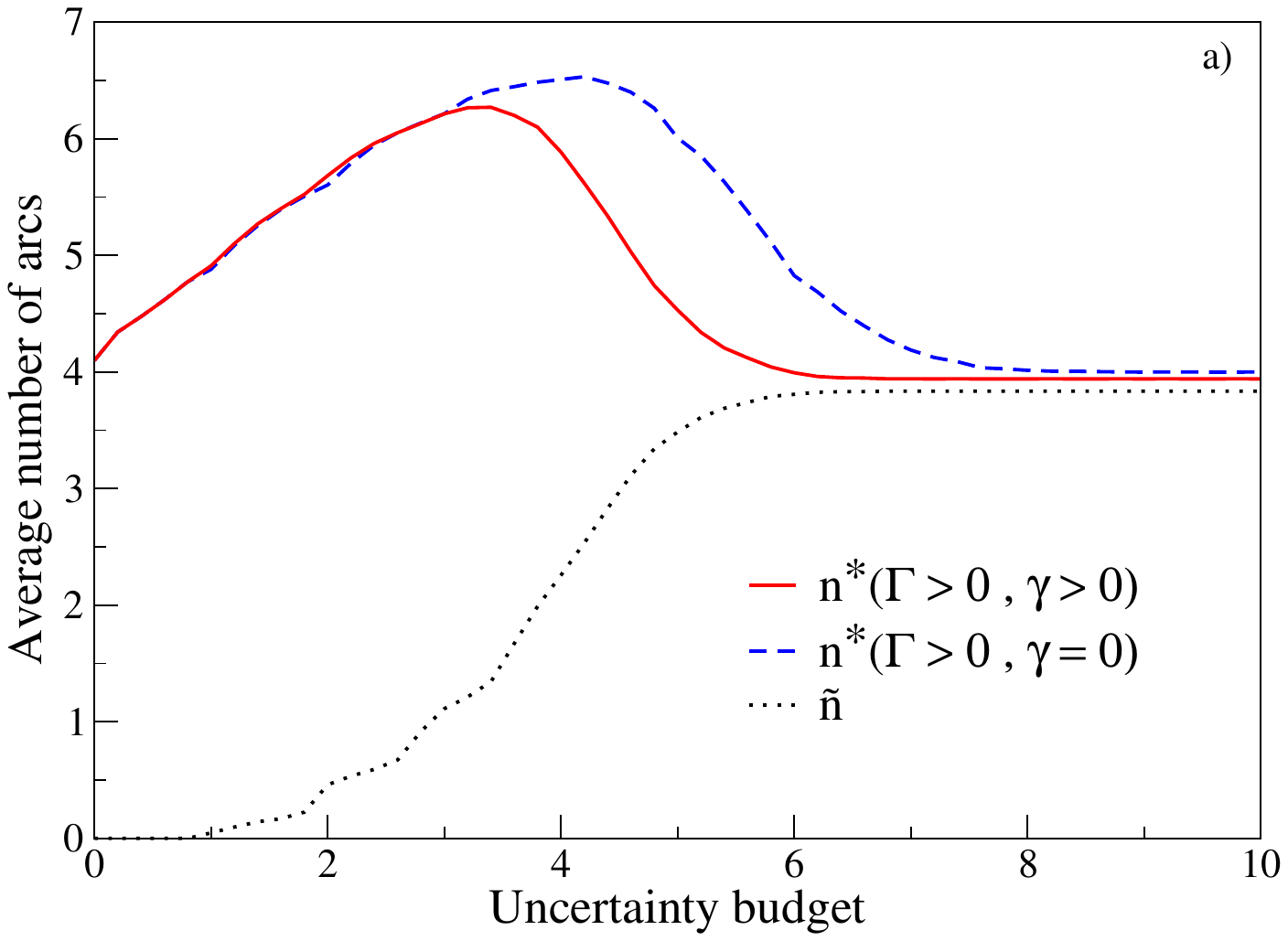}
	   \end{minipage}
	\hspace{-5mm}
  	   \begin{minipage}[t]{75mm}
		\includegraphics[width=75mm]{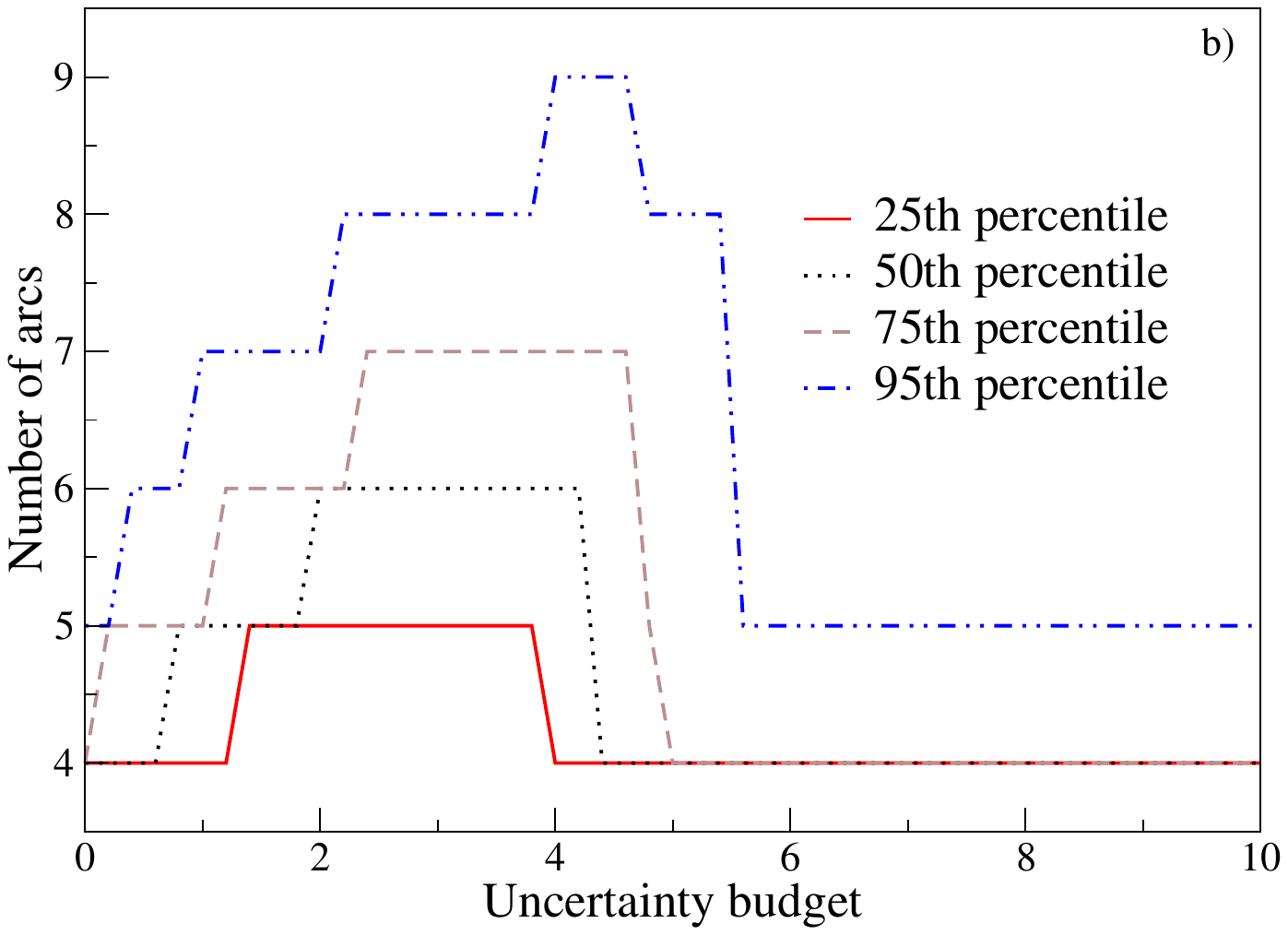}
 	\end{minipage}
\caption{The dependence on budget of uncertainty $\Gamma$ for: a) average number of arcs and b) their distribution. The graph consists of \(|\mathcal{V}|=30\) nodes and uncertainty reduction is permitted.}
\label{figr:46}
\end{figure}

Figure~\ref{figr:46}a displays how the average $n^*$ changes with \(\Gamma\) for the different settings.
Note that the values of uncertainty are relative to the nominal arc length. 
This provides an upper bound on the maximum objective value, i.e., when every arc in the shortest path (contributing to \(n^{*}\)) is affected by the uncertainty.
At $\Gamma=0$, we observe $n^*(\Gamma=0,\gamma=0)$, and $\tilde{n}=0$.
As $\Gamma$ increases, it turns beneficial to choose more but shorter arcs, hence, the average \(n^{*}(\Gamma>0,\gamma=0)\) initially increases and reaches a maximum at $\Gamma\approx n^*(\Gamma=0,\gamma=0)$.
As $\Gamma$ grows even further, the standard robust solution $n^*(\Gamma>0,\gamma=0)$ decreases and plateaus at the same level as $n^*(\Gamma=0,\gamma=0)$.
When $\gamma>0$, we observe that an increasing $\Gamma\ge0$ permits more uncertain arc lengths to be reduced ($\tilde{n}\ge0$) to a maximum of $\tilde{n}\le n^*(\Gamma=0,\gamma=0)$. 
Since some of the arc uncertainty can be reduced, the peak of $n^*(\Gamma>0,\gamma>0)$ occurs at a lower budget than when no reduction is allowed, as seen in Figure~\ref{figr:46}a.
Note that for small \(\Gamma\), in order to cope with uncertainty, the optimal solution minimizes the length of each individual arc so that the impact of the uncertainty is minimized. 

To further support this observation, Figure~\ref{figr:46}b displays the distribution of the number of arcs using different percentiles of $n^*(\Gamma>0,\gamma>0)$ (corresponding to Figure~\ref{figr:46}a). 
Here, we observe that as \(\Gamma\) increases, the distribution of $n^*(\Gamma>0,\gamma>0)$ skews towards larger number of arcs (the gaps between the percentiles increase).
This means that the optimal solution becomes more diversified.
Specifically, the model selects a path consisting of some certain and some uncertain arcs, with a subset of the latter experiencing uncertainty reduction.
This continues until the saturation point (here $\Gamma\approx4$) because
beyond a certain budget, diversification of paths becomes redundant.
At this point, the shortest path is chosen exclusively amongst uncertain arcs, almost all of which experience uncertainty reduction (since $\Gamma> n^*(\Gamma=0,{\gamma=0})$).

\noindent\emph{Experiment 6: Comparison to SO. $\quad$}
This experiment evaluates the average and worst case performance of the robust DDU solutions and compares them to a similar SO problem. 
The SO formulation is given by
\begin{eqnarray*}
  \min_{\mb{x}}&&\; \sum_{(i,j) \in \mathcal{A}}c_{ij}x_{ij} + \mathbb{E}_{\mathbb{P}(\mb{x})}\left[\sum_{({i,j}) \in \mathcal{A}} d_{ij}(\sbf{\xi})y_{ij}\right]\\
  \text{s.t.}&&\; \mb{y} \in Y\\
          &&\; \mb{x} \in \{0,1\}^{|\mathcal{A}|},
\end{eqnarray*}
with the uncertainty set
\[\bs{\xi} \in \mU^{SSP}(\mb{x}) = \times_{{i,j} \in \mathcal{A}} [0 , 1 - \gamma x_{ij}].\]
The distribution \(\mathbb{P}(\mb{x})\) is the uniform distribution over the support \(\mU^{SSP}(\mb{x})\).
The average performance is evaluated by randomly generating the uncertain component \(\xi_{ij}\) (from \([0,1]\) for unreduced arcs and \([0,1-\gamma x_{ij}\)] for reduced arcs) and implementing the existing robust and stochastic solutions for these randomly generated arc costs.
The following solutions are evaluated:
\begin{enumerate*}[label=(\roman*)]
	\item RO: Robust solution for $\gamma=0$. 
	\item RO-DDU: Robust solution for $\gamma>0$.
	\item SO: Stochastic solution for $\gamma=0$. 
	\item SO-DDU: Stochastic solution for $\gamma>0$. 
\end{enumerate*}
The suffix of the average performances is ``-A'' and of the worst case performances ``-W.''
\begin{figure}[h!]
\centering
\includegraphics[width=75mm]{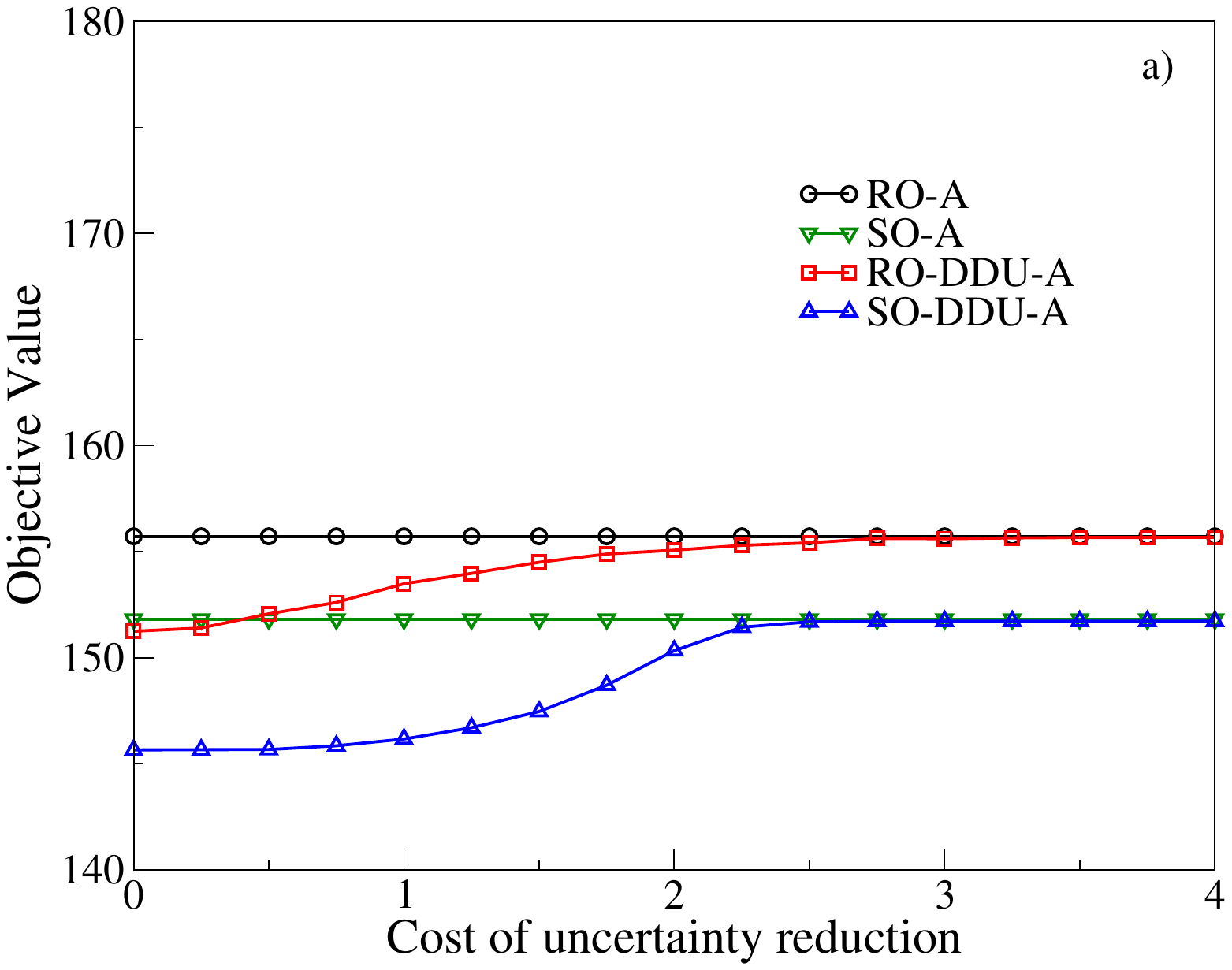}
\hspace{-3mm}
\includegraphics[width=75mm]{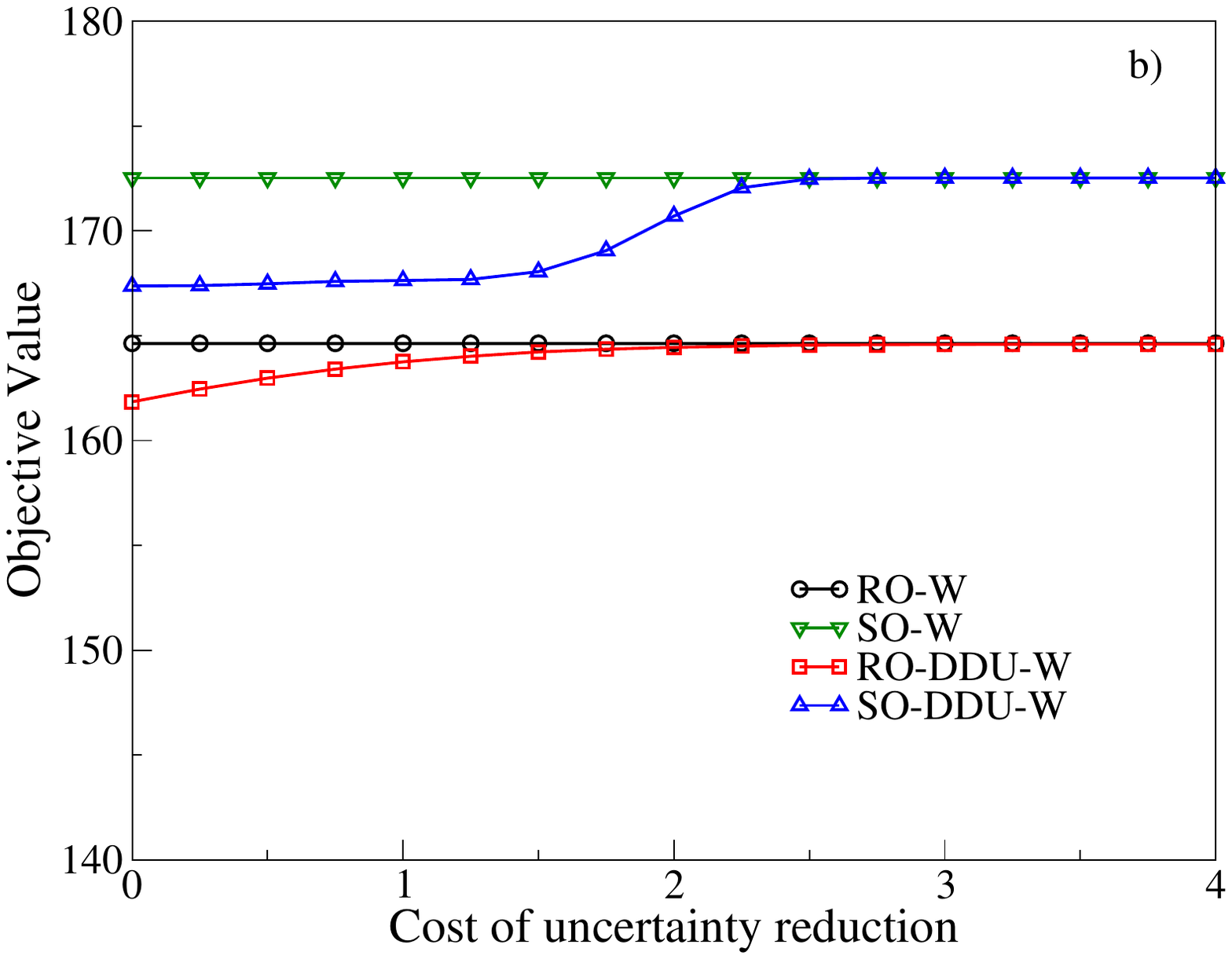}
\caption{Comparison of RO and SO formulations: a) average and b) worst-case objective value.}
\label{figr:SG_AaW}
\end{figure}

Figure~\ref{figr:SG_AaW}a shows that the average objective of SO is less than the average RO objective. 
This is because RO optimizes the worst-case instead of the average performance as in SO\@.
However, analogously in Figure~\ref{figr:SG_AaW}b, RO-W is significantly less than SO-W.
The same applies to the \dd counterparts for both cases.
As can be expected, the objective values increase with \(c\) until it is no longer beneficial to reduce the uncertainty, i.e., the objective value of the RO-DDU solution increases until it matches that of the RO solution.
The same holds true for the SO-DDU and SO solutions.

In summary, the \(\overline{\bs{\Pi}}\)-formulation and the Modified \mbox{Big-M} formulation perform considerably better than the standard \mbox{Big-M} formulation and their benefits increase with graph size. 
The worst-case cost for the shortest path can be improved by proactively reducing the uncertainty on a subset of arcs. 
As the budget of uncertainty grows, these benefits improve but plateau beyond a certain level. 
At the same time, the cost of reduction curbs these benefits. 
The RO-DDU problem performs better than SO-DDU for the worst-case scenario.
As expected, this benefit comes at the price of the average cost. 
This numerical study provides an overview of the impact of different formulations, probes various model parameters, and highlights the power of the proactive uncertainty control for both the worst-case and average performance.

%%%%%%%%%%%%%%%%%%%%%%%%%%%%%%%%%%%%%%%%%%%%%%%%%%%%%%%%%%%%%%%%%%%%%%%%%%%%%%%

\section{Concluding Remarks}
\label{sec:conclusion}
In this paper, we present a novel optimization approach for solving problems with \dd uncertainties.
We show that for general polyhedral sets, such problems are, even in basic cases, NP-complete.
To alleviate this, we introduce a class of uncertainty sets whose upper bounds are affected by decisions.
They enable more realistic modeling of a broad range of applications and extend RO beyond the currently used exogenous sets.
We provide reformulations that have considerably fewer constraints compared to standard linearization techniques, allowing for faster computations.
Our approach should be viewed as one option among many to model decision dependence while maintaining computational advantages.
The induced convexity of the sub-problem in the proposed reformulation reveals a path forward to use advanced cut generating algorithms.
We believe that finding new and appropriate conditions on sets will further improve the quality of the reformulations.

In addition, this work provides an alternative way of addressing one of the criticisms of RO approaches, namely overly conservative solutions.
The description via \dd sets enables mitigation of this issue by exercising proactive control on uncertainties.
This setting offers an immediate way to manage the tradeoff between conservatism and optimality.
Finally, novel cutting plane methods have instrumentally enhanced solution times and we envision \dd sets to solidify the tradeoff between computation and optimality by inducing beneficial cuts.

%\appendix

\section*{Acknowledgments}
We are grateful to David Morton for insightful comments.

\bibliographystyle{abbrvnat-onn}
\bibliography{DDUnewpaper_short}

\end{document}